\documentclass[a4paper, reqno]{amsart}
\usepackage{amsmath, amssymb, eucal, amscd, amstext, enumerate, mathrsfs, yfonts, amsfonts, verbatim, tikz}

\newtheorem{thm}{Theorem}[section]
\newtheorem{lem}[thm]{Lemma}
\newtheorem{eg}[thm]{Example}
\newtheorem{prop}[thm]{Proposition}
\newtheorem{prop-defn}[thm]{Proposition and Definition}
\newtheorem{cor}[thm]{Corollary}

\newtheorem{rem}[thm]{Remark}

\newenvironment{prf}{{\noindent \textbf{Proof:}\ }}{\hfill $\Box$\\ \smallskip}

\numberwithin{equation}{section}

\newcommand{\ti}{\tilde}
\newcommand{\smnoind}{\smallskip\noindent}
\newcommand{\CL}{\mathcal{L}}
\newcommand{\la}{\langle}
\newcommand{\ra}{\rangle}
\newcommand{\CH}{\mathcal{H}}
\newcommand{\CK}{\mathcal{K}}

\newcommand{\CC}{\mathcal{C}}
\newcommand{\CT}{\mathcal{T}}
\newcommand{\BR}{\mathbb{R}}
\newcommand{\BP}{\mathbb{P}}
\newcommand{\RP}{\mathbb{R}_+}
\newcommand{\BN}{\mathbb{N}}

\newcommand{\KM}{\mathfrak{M}}
\newcommand{\KG}{\mathfrak{G}}

\newcommand{\kk}{\mathfrak{k}}
\newcommand{\um}{\mathbf{u}}
\newcommand{\id}{\mathrm{id}}
\newcommand{\BC}{\mathbb{C}}
\newcommand{\KI}{\mathfrak{I}}

\newcommand{\Ber}{\mathbf{Ber}}

\newcommand{\KJ}{\mathfrak{J}}
\newcommand{\BZ}{\mathbb{Z}}
\newcommand{\BA}{\mathbb{A}}
\newcommand{\TBK}{{\BA_{\BZ_1}^1}}
\newcommand{\TBKM}{\BA_{\BZ_1}^{1,\min}}
\newcommand{\AZ}{{\BA_{\BZ_1}^1}}
\newcommand{\AZU}{{\BA_{\BZ_0}^{1}}}
\newcommand{\BK}{\mathbb{K}}
\newcommand{\BQ}{\mathbb{Q}}
\newcommand{\Aut}{\mathrm{Aut}}
\newcommand{\bt}{\mathbf{t}}
\newcommand{\KPI}{{\BZ[\bt]}}
\newcommand{\bp}{\mathbf{p}}
\newcommand{\bq}{\mathbf{q}}
\newcommand{\kp}{\mathfrak{P}}
\newcommand{\E}{1}

\begin{document}

\title{Berkovich spectra of elements in Banach Rings\footnote{This work is supported by the National Natural Science Foundation of China (11071126 and 11471168)}
}

\author{Chi-Wai Leung and Chi-Keung Ng}

\address[Chi-Wai Leung]{Department of Mathematics, The Chinese University of Hong Kong, Hong Kong.}
\email{cwleung@math.cuhk.edu.hk}

\address[Chi-Keung Ng]{Chern Institute of Mathematics and LPMC, Nankai University, Tianjin 300071, China.}
\email{ckng@nankai.edu.cn}

\date{\today}

\begin{abstract}
Adapting the notion of the spectrum $\Sigma_a$ for an element $a$ in an ultrametric Banach algebra over a complete valuation field (as defined by Berkovich), we introduce and briefly study the Berkovich spectrum $\sigma^\Ber_R(u)$ of an element $u$ in a Banach ring $R$.
This spectrum is a compact subset of the affine analytic space $\AZ$ over $\BZ_1$ (the ring $\BZ$ equipped with the Euclidean norm), and the later can be identified with the ``equivalence classes'' of all elements in all complete valuation fields.
If $R$ is generated by $u$ as a unital Banach ring, then $\sigma^\Ber_R(u)$ coincides with the spectrum $\KM(R)$ of $R$.
If $R$ is a unital complex Banach algebra, then $\sigma^\Ber_R(u)$ is the ``folding up'' of the usual spectrum $\sigma_B(u)$ alone the real axis.

For a non-Archimedean complete valuation field $\kk$ and an infinite dimensional ultrametric $\kk$-Banach space $E$ with an orthogonal base, if $u\in \CL(E)$ is a completely continuous operator, we show that many different ways to define the spectrum of $u$ give the same compact set
$\sigma^\Ber_{\CL(E)}(u)$.
As an application, we give a lower bound for the valuations of the zeros of the Fredholm determinant $\det(1-\bt u)$ (as defined by Serre) in complete valuation field extensions of $\kk$.
Using this, we give a concrete example of a completely continuous operator whose Fredholm determinant does not have any zero in any complete valuation field extension of $\kk$.

On our way, we also give a  complete description of the topological subspace $\TBKM\subseteq \AZ$ consisting of homeomorphic images of all ``minimal'' complete valuation fields.

\end{abstract}
\maketitle

\section{Introduction}

In Chapter 1 of \cite{Berk90}, Berkovich defined the spectrum $\KM(T)$ for a commutative unital Banach ring $T$.
He used it as one of the tools for his non-Archimedean geometry.
He also defined in \cite[Chapter 7]{Berk90}, the spectrum $\Sigma_f$ of an element $f$ in an ultrametric Banach algebra $A$.
In this article, we define and study the spectrum $\sigma^\Ber_R(a)$ of an element $a$ in a general unital Banach ring $R$.

\medskip

As in \cite{Berk90}, a canonical way to define such a spectrum is to consider those elements $t$ in a complete valuation field $K$ such that $a\otimes 1 - 1\otimes t$ is not invertible in the (projective) Banach ring tensor product $R\hat \otimes_\BZ K$.
There are two issues that one needs to pay attention at.
The first one is that $R\hat \otimes_\BZ K$ could be zero.
This issue can easily be resolved by adding the assumption that it is non-zero, in the definition of the spectrum.
The second issue is that if $K$ is a valuation subfield of another complete valuation field $L$, we want to identify $t\in K$ with $t\in L$.

\medskip

This second issue gives rise to the following equivalent relation:
$(K_1,t_1)\sim (K_2,t_2)$ if there exists a complete valuation field $K$ and an element $t\in K$ such that $K$ can be identified with valuation subfields of both $K_1$ and $K_2$ under which $t$ equals $t_1$ and $t_2$, respectively.
In fact, the equivalence classes of the disjoint union of all elements in all complete valuation fields (strictly speaking, this union is not a set) is a set with a natural topology, under which it becomes a second countable locally compact Hausdorff space.

\medskip

After defining and studying the basic property of the spectra of elements in unital Banach rings as well as their counterparts for Banach algebras, we will use it to tackle the question of whether the Fredholm determinant of a completely continuous operator has a zero in a complete valuation field extension.
This question is another motivation of this study.

\medskip

The paper is organized as follows.
Since (projective) tensor products of Banach rings play an important role in our study and this notion is not well-documented in the literature, we will give a brief account of it in Section 2.
Moreover, we will use the tensor product construction to show that a unital Banach ring is a Banach subring of a unital Banach algebra if and only if it satisfies a canonical ``regularity condition'' (Proposition \ref{prop:reg-ring-ten-prod}).
We will also recall in that section the definition and basic properties of $\KM(T)$.

\medskip

In Section 3, we will study the space of equivalence classes of all elements in all complete valuation fields.
More precisely, it is established in Proposition \ref{prop:top-K-alpha} that if $\alpha$ is a fixed infinite cardinal and $\CK$ is the set of all complete valuation fields with cardinalities less than $2^\alpha$, then there is a bijective correspondence $\Theta$ from the set $\ti \BK$ of equivalence classes of the disjoint union of fields in $\CK$ (under the equivalence relation as defined above) to the set $\AZ$ of non-zero multiplicative semi-norms on the ring $\KPI$ of polynomials in one variable with integral coefficients (note that the later is independent of $\alpha$).
It is well-known that the Hausdorff topology on $\AZ$ given by pointwise converges on elements in $\KPI$ is second countable and locally compact.
Moreover, there is a canonical ``valuation'' on $\AZ$ given by $\lambda \mapsto \lambda(\bt)$ under which $\Theta$ preserves the valuations induced from those valuation fields.
%It is found that subsets of $\AZ$ that are closed and bounded (with respects to this valuation)  are compact.
We will identify $\ti \BK$ with $\AZ$ directly.

\medskip

Suppose that $\mu$ is the canonical continuous map from the disjoint union of elements in $\CK$ to $\TBK$.
The restriction of $\mu$ on $\BC$ can be regarded as the ``folding up'' of the complex
plane alone the real axis (Example \ref{eg:non-inj}(b)).
If $K$ is non-Archimedean, its image under $\mu$ is inside the subset, $\AZU$, consisting of ultrametric semi-norms.

\medskip

The restrictions of $\mu$ on minimal fields in $\CK$ are homeomorphisms onto their images and all such images are disjoint from one another (see Proposition \ref{prop:top-K-alpha} and Lemma \ref{lem:top-K-min}).
In Theorem \ref{thm:descr-K-min}, we will give a complete description of the topological subspace $\TBKM$ of the unions of all such images.

\medskip

In Section 4, we will give the definition and some study of the Berkovich spectrum $\sigma^\Ber_R(a)$ as well as its ultrametric counterpart $\sigma^\um_R(a)$.
Using the idea in the proof of the compactness of $\Sigma_f$ in \cite{Berk90}, we will show in Theorem \ref{thm:closed-ring-sp} that both $\sigma^\Ber_R(a)$ and $\sigma^\um_R(a)$ are compact subsets of $\TBK$.
We will also show that $\sigma^\Ber_R(a)$ is non-empty if there is a contractive additive map from $R$ to a complete valuation field, while $\sigma^\um_R(a)$ is non-empty if and only if there exists a contractive additive map from $R$ to a non-Archimedean complete valuation field.
In particular, if $R$ is any unital ring equipped with the trivial valuation, then $\sigma^\um_R(a)$ is always a non-empty subset of $\AZU$ (Theorem \ref{thm:non-empty:ring-triv-val}(b)).
If $R$ is generated by $a$ as a unital Banach ring, then Proposition \ref{prop:Omega->sigma}(b) tells us that $\sigma^\Ber_R(a)$ coincides with $\KM(R)$.
However, unlike the cases of complex Banach algebras and ultrametric Banach algebras, if $S$ is the unital Banach subring of $R$ generated by $a$, it is possible that $\partial \sigma_S^\Ber(a) \nsubseteq \sigma_R^\Ber(a)$ even when $R$ is commutative (see Example \ref{eg:spec-diff-rings}). 

\medskip

Furthermore, we consider the case when $R$ is a Banach algebra over a complete valuation field $\kk$.
We define closed subsets $\sigma^\Ber_{R,\kk}(a)$ and $\sigma^\um_{R,\kk}(a)$ of $\sigma^\Ber_R(a)$ and $\sigma^\um_R(a)$, respectively, that take into the account of the scalar multiplication on $R$.
It is stated in Proposition \ref{prop:cpt-non-empty-alg-spec} that $\sigma^\Ber_{R,\kk}(a)$ (respectively, $\sigma^\um_{R,\kk}(a)$) is non-empty if there exists a non-zero contractive
additive map from $R$ to $\kk$ (respectively, and $\kk$ is ultrametric).
Therefore, if $R$ is the unitalization of another $\kk$-Banach algebra, then $\sigma^\Ber_{R,\kk}(a)\neq \emptyset$ (Proposition \ref{prop:cpt-non-empty-alg-spec}(f)).
Moreover, when $R$ is ultrametric, $\sigma^\um_{R,\kk}(a)$ is the canonical image of $\Sigma_a$ in $\AZU$ (Lemma \ref{lem:rel-bet-Sigma-and-Ber}).
In the case when $R$ is a complex Banach algebra, one has $\sigma^\Ber_R(a) = \sigma^\Ber_{R,\BC}(a)$ and it is the ``folding up'' of the ordinary spectrum $\sigma_R^\BC(a)$ alone the real axis (see Proposition \ref{prop:spec-elem-in-real-alg} and Example \ref{eg:non-inj}(b)).

\medskip

In Section 5, we will consider the case when $\kk$ is non-Archimedean, $R$ is the $\kk$-Banach algebra $\CL(E)$ of bounded linear maps on an infinite dimensional ultrametric $\kk$-Banach space $E$ with an orthogonal base, and $a\in R$ is a completely continuous operator.
In this case, it is shown in Proposition \ref{prop:inv-in-A} that $\sigma^\um_{R,\kk}(a) = \sigma^\Ber_R(a)$ and they also coincide with many different notions of the spectrum of $u$.
Using these, we obtain a relationship between ``non-zero'' elements of $\sigma^\Ber_R(a)$ and zeros of the Fredholm determinant $\det(1-\bt a)$ of $a$ (Theorem \ref{thm:zero-Fred-det}), and show that if $\lambda$ is a zero of $\det(1-\bt a)$ in a complete valuation field extension $K$ of $\kk$, then $|\lambda|_K \geq \lim_n \|a^n\|^{-1/n}$ (Corollary \ref{cor:lower-bdd-zero}). 
We will then give a concrete example of a completely continuous operator whose Fredholm determinant has no zero in any complete valuation field extension of $\kk$ (Example \ref{eg:Fredh-det-no-zero}).

\medskip

\section{Notation and Preliminary on Banach rings}

In this article, we denote $\BN_0:=\BN\cup \{0\}$.
A \emph{semi-norm} on an additive (i.e.\ abelian) group $X$ is a subadditive function $\|\cdot\|: X \to \RP$ satisfying $\|-x\| = \|x\|$ ($x\in X$).
We set $\ker \|\cdot\| := \{x\in X: \|x\| = 0\}$.
A semi-norm $\|\cdot\|$ is said to be \emph{ultrametric} if $\|x+y\| \leq \max\{\|x\|,\|y\|\}$.
On the other hand, it is called a \emph{norm} if $\ker\|\cdot\| = \{0\}$.
For any additive groups $X$ and $Y$ with fixed semi-norms, an additive map $\varphi:X\to Y$ is said to be \emph{contractive} if $\|\varphi(x)\| \leq \|x\|$ ($x\in X$).
If $\|\cdot\|$ is a norm on $X$ and $X$ is complete under the metric defined by $d(x,y):=\|x-y\|$, then $X$ is called a \emph{Banach additive group}.
One can always ``complete'' a normed additive group to obtain a Banach additive group.

\medskip

Suppose that $R$ is a unital ring.
We denote by $R[\bt]$ (respectively, $R[[\bt]]$) the ring of all polynomials (respectively, formal power series) in one variable with coefficients in $R$.
Moreover, we consider $R_0$ to be the ring $R$ equipped with the trivial norm $\|\cdot\|_0$ (namely, $\|x\|_0 := 1$ for every $x\in R\setminus \{0\}$). 

\medskip

If $R$ is also a Banach additive group with the same addition such that the norm is submultiplicative, then we call $R$ a \emph{Banach ring}.
We use the term complete valuation fields, non-Archimedean valuation fields, Banach modules and Banach spaces in their usual senses.
If $R$ is a Banach ring with an identity $1$, we say that it is a \emph{unital Banach ring} if $\|1\| = 1$.

\medskip

For any unital Banach rings $R$ and $S$, we denote by $\CC(R;S)$ the set of all contractive unital ring homomorphisms from $R$ to $S$, and by $\Aut(R)$ the set of bijective elements in $\CC(R;R)$.
If $R$ and $T$ are unital Banach rings with $T$ being commutative, we say that $R$ is a \emph{unital Banach $T$-algebra} if there exists $\varphi\in \CC(T;R)$; in this case, we set $s\cdot x := \varphi(s)x$ ($s\in T;x\in R$).
Notice that a unital Banach ring is a unital Banach $\BZ_1$-algebra, where $\BZ_1$ is the ring $\BZ$ equipped with the Euclidean norm.

\medskip

Suppose that $K$ and $L$ are complete valuation fields.
We say that $L$ is an \emph{extension} of $K$ if $K$ is isometrically isomorphic to a norm closed subfield of $L$.
Moreover, $K$ is said to be \emph{minimal} if $K$ is the only complete valuation subfield contained in $K$.
We denote by $\CK^{\min}$ the set of all minimal complete valuation fields.
By the Ostrowski's theorem, elements in $\CK^{\min}$ are algebraically isomorphic to either $\BQ$, $\BR$, $\BZ/p\BZ$ or $\BQ_p$, for some prime number $p$.
%Moreover, in the case of $\BQ$ or $\BZ/p\BZ$, only the trivial valuation is considered, while in the case of $\BR$ (respectively, $\BQ_p$), all valuations that are equivalent to the Euclidean valuation (respectively, the $p$-adic valuation) are considered, and they are different elements in $\CK^{\min}$.
%
%\medskip
%
More precisely, for any $\upsilon\in (0,1]$, wet set $\BR_\upsilon$ to be the field $\BR$ equipped with the valuation $t\mapsto |t|_\E^\upsilon$ (where $|\cdot|_\E$ is the Euclidean norm on $\BR$).
Similarly, for any $p\in \BP$ and $\omega\in (0,\infty)$, we set $\BQ_p^\omega$ to be the field $\BQ_p$ equipped with the valuation $t\mapsto |t|_p^\omega$ (where $|\cdot|_p$ is the $p$-adic norm on $\BQ_p$).
%In particular, $\BR_1$ and $\BQ_p^1$ are the fields $\BR$ and $\BQ_p$ equipped with the usual Euclidean norm and the $p$-adic norm, respectively.
We will also use $\BQ_0$ and $\BZ(p)$ to denote the fields $\BQ$ and $\BZ/p\BZ$, both equipped with the trivial valuation.
If we set,
$\CK^{\min}_\BR := \{\BR_\upsilon: \upsilon\in (0,1]\}$ and $\CK^{\min}_{\BQ_p} := \{\BQ_p^\omega: \omega\in (0,\infty)\}$,
then
\begin{equation}\label{eqt:decomp-min-val-field}
\CK^{\min} = \{\BQ_0\} \cup \CK^{\min}_\BR \cup \{\BZ(p):p\in \BP\}\cup \bigcup_{p\in \BP} \CK^{\min}_{\BQ_p},
\end{equation}
where $\BP$ is the set of all prime numbers.

\medskip

On the other hand, we say that a complete valuation field $K$ is \emph{generated by a subset $E$} if the smallest closed subfield of $K$ containing $E$ is $K$ itself.
If it happens that $E = \{r\}$, we say that $K$ is \emph{singly generated} and that $r$ is a \emph{generator} for $K$.
Clearly, any minimal field is singly generated but the converse is not true (e.g.\ $\mathrm{i}$ is a generator for $\BC$).

\medskip

The following are some easy facts, whose proof are left to the readers.

\medskip

\begin{lem}\label{lem:k-lin}
Suppose that $K$ and $L$ are complete valuation fields.
Let $R$ be a unital Banach ring and $S$ be a unital ring. 

\smnoind
(a) Any $\varphi\in \CC(K;R)$ is an isometry.

\smnoind
(b) If $K\in\CK^{\min}$, then $\CC(K;K) = \{\id\}$.

\smnoind
(c) If $\gamma$ is a submultiplicative semi-norm on $S$ such that there is $\kappa\in \RP$ with  
$\gamma(s^n) = \gamma(s)^n$ ($s\in S;n\in \BN$) and $\gamma(m\cdot 1)\leq \kappa$ ($m\in \BZ$), then for any $x,y\in S$ satisfying $xy=yx$, one has 
$\gamma(x+y) \leq \max\{\gamma(x), \gamma(y)\}$.

\smnoind
(d) Let $\KG(R)$ be the set of invertible elements in $R$. 
If $x\in R$ satisfying $\|1-x\| < 1$, then $x\in \KG(R)$.
Consequently, $\KG(R)$ is an open subset of $R$.
\end{lem}
%\begin{comment}
[This proof will not appear in the published version.] 

\begin{prf}
(a) For any $s\in K\setminus \{0\}$, we have 
$$\|\varphi(s)\|^{-1}\ \leq\ \|\varphi(s^{-1})\|\ \leq\ |s^{-1}|\ =\ |s|^{-1}\ \leq\ \|\varphi(s)\|^{-1},$$
which shows that $\varphi$ is an isometry.

\smnoind
(b) By part (a), any $\varphi\in \CC(K;K)$ will have a closed image, which is closed subfield of $K$ (which is $K$ itself). 
Moreover, since $\varphi(m) = m$ ($m\in \BZ$), the restriction on $\varphi$ on the subfield generated by the image of $\BZ$ in $K$ (which is $K$ itself) is the identity map. 

\smnoind
(c) Without loss of generality, we may assume that $\kappa \geq 1$.
If $M := \max\{\gamma(x), \gamma(y)\}$, then for any $n\in \BN$, one has 
$\gamma(x+y)^n  = \gamma(\sum_{k=0}^n C^n_k x^k y^{n-k}) \leq \sum_{k=0}^n \gamma(C^n_k \cdot 1) M^n \leq (n+1) \kappa M^n$. 
This gives $\gamma(x+y) \leq M$ as required. 

\smnoind
(d) The first statement follows from the argument for a similar statement concerning Banach algebra. 
The second statement follows from $\|x^{-1}y -1 \|, \|y x^{-1} - 1\| \leq \|x^{-1}\| \|y-x\|$. 
\end{prf}
%\end{comment}

\medskip

Notice that there are exactly two elements in $\CC(\BC;\BC)$ (one of them is given by $s \mapsto \bar s$) and $\BC$ can be regarded as an extensions of $\BC$ in two different ways.

\medskip

Suppose that $T$ is a commutative unital Banach ring.
Let $X$ and $Y$ be Banach $T$-modules.
We define a semi-norm $\|\cdot\|_\wedge$ and an ultrametric semi-norm $\|\cdot\|_\wedge^\um$ on the algebraic tensor product $X\otimes_T Y$ over $T$ by
$$\|z\|_\wedge := \inf \Big\{\sum_{k=1}^n \|a_k\|\|b_k\|: n\in \BN; a_1,...,a_n\in X; b_1,...,b_n\in Y; z = \sum_{k=1}^n a_k\otimes b_k\Big\}$$
and
$$\|z\|_\wedge^\um := \inf \Big\{\max_{k=1,...,n} \|a_k\|\|b_k\|: n\in \BN; a_1,...,a_n\in X; b_1,...,b_n\in Y; z = \sum_{k=1}^n a_k\otimes b_k\Big\}.$$
We denote by $X\hat\otimes_T Y$ and $X\hat\otimes_T^\um Y$
the completions of $X\otimes_T Y / \ker\|\cdot\|_\wedge$ and $X\otimes_T Y / \ker\|\cdot\|_\wedge^\um$, respectively.
In this case, both $X\hat\otimes_T Y$ and $X\hat\otimes^\um_T Y$ are Banach $T$-modules.
By abuse of notation, we identify $x\otimes y$ with its images in $X\otimes_T Y$,  $X\hat\otimes_T Y$ and $X\hat\otimes_T^\um Y$.

\medskip

The argument of the following result is standard and is left to readers.

\medskip

\begin{lem}\label{lem:univ-prop-proj-ten}
Let $T$ a commutative unital Banach ring, and let $X$, $Y$ and $Z$ be unital Banach $T$-modules.

\smnoind
(a) $X\hat\otimes_T Y$ is a Banach $T$-module and $X\hat\otimes_T^\um Y$ is an ultrametric Banach $T$-module.

\smnoind
(b) Suppose that $\varphi: X\times Y\to Z$ is a map such that
$\varphi(rx,y) = \varphi(x,ry) = r\varphi(x,y)$ and $\|\varphi(x,y)\| \leq \|x\|\|y\|$ ($x\in X;y\in Y;r\in T$).
There is a unique contractive $T$-module map $\hat \varphi: X\hat \otimes_R Y\to Z$ with $\hat\varphi(x\otimes y) = \varphi(x,y)$.
If, in addition, the norm on $Z$ is ultrametric, there is a unique contractive $T$-module map $\hat \varphi: X\hat \otimes^\um_R Y\to Z$ with $\hat\varphi^\um(x\otimes y) = \varphi(x,y)$.
These two universal properties completely characterize the Banach $T$-modules $X\hat \otimes_R Y$ and $X\hat \otimes^\um_R Y$ (together with the canonical embeddings), respectively.

\smnoind
(c) Suppose that $X$, $Y$ and $Z$ are unital Banach $T$-algebras.
Then $X\hat\otimes_T Y$ (respectively, $X\hat\otimes_T^\um Y$) is either zero or a unital (respectively, ultrametric) Banach $T$-algebra.
If the map $\varphi$ in part (b) also satisfies $\varphi(1,1) = 1$ and  $\varphi(ax,by) = \varphi(a,b)\varphi(x,y)$ ($a,x\in X;b,y\in Y$), then $\hat\varphi$ (respectively, $\hat\varphi^\um$) is a unital $T$-algebra homomorphism.
Furthermore, this universal property concerning contractive $T$-algebra homomorphisms will characterise the unital Banach $T$-algebra $X\hat\otimes_T Y$ (respectively, $X\hat\otimes_T^\um Y$).
\end{lem}

\medskip

\begin{lem}\label{lem:proj-ten}
Suppose that $R$ and $S$ are unital Banach rings. 

\smnoind
(a) Let $K$ be a complete valuation field. 
If $R\hat \otimes_\BZ K \neq (0)$ (respectively, $R\hat \otimes_\BZ^\um K \neq (0)$), then  it is a unital $K$-Banach algebra.

\smnoind
(b) If $K\in \CK^{\min}$ and both $R$ and $S$ are unital $K$-Banach algebras, then $R\hat\otimes_\BZ S \cong R\hat\otimes_K S$ and $R\hat\otimes^\um_\BZ S \cong R\hat\otimes^\um_K S$ canonically.
\end{lem}
\begin{prf}
(a) This follows from Lemmas \ref{lem:k-lin}(a) and \ref{lem:univ-prop-proj-ten}(c).

\smnoind
(b) If $K$ is of characteristic $p\in \BP$ (i.e., $K = \BZ(p)$), then 
$R\otimes_\BZ S = R\otimes_K S$ and the corresponding projective tensor norms on these two algebras are the same.
Suppose that $K$ is of characteristic zero.
Then $R\otimes_\BZ S = R\otimes_\BQ S$.
Moreover, since any $s\in K$ can be approximated by a sequence $\{s_n\}_{n\in \BN}$ in $\BQ$ and
$$\|as\otimes b - a\otimes sb\|_{R\hat\otimes_\BZ S} \leq \|a(s-s_n)\otimes b\|_{R\hat\otimes_\BZ S} + \|a\otimes (s_n-s)b\|_{R\hat\otimes_\BZ S}
\quad (a\in R; b\in S),$$
we know that $as\otimes b = a\otimes sb$ in $R\hat\otimes_\BZ S$.
Thus, Lemma \ref{lem:univ-prop-proj-ten}(c) produces a map in $\CC(R\hat\otimes_K S;R\hat\otimes_\BZ S)$ that respects simple tensors, which means that the canonical map from $R\hat\otimes_\BZ S$ to $R\hat\otimes_K S$ is a bijective isometry.
The argument for $R\hat\otimes^\um_\BZ S \cong R\hat\otimes^\um_K S$ is similar.
\end{prf}

\medskip

\begin{lem}\label{lem:non-zero-ten-norm-sp}
Let $K$ be a non-Archimedean complete valuation field and $K^0$ be the smallest closed subfield contained in $K$.
Let $E$ and $F$ be $K$-Banach spaces.

\smnoind
(a) If $F$ is ultrametric, then $\|\cdot\|_\wedge$ is a norm on $E\otimes_K F$.

\smnoind
(b) The following statements are equivalent.
\begin{enumerate}
\item There is a non-zero contractive additive map $\phi$ from $E$ to $K$.
%\item There is an ultrametric $K$-Banach space $E_0$ and a non-zero contractive $K$-linear map from $E$ to $E_0$.
\item $E\hat\otimes_K^\um K \neq (0)$.
\item $E\hat\otimes_\BZ^\um K^0 \neq (0)$.
\end{enumerate}

\smnoind
(c) $\|\cdot\|_\wedge^\um$ is a norm on $E\otimes_K K$ if and only if for any $x\in E$, one can find a contractive additive map $\phi: E \to K$ with $\phi(x)\neq 0$.
\end{lem}
\begin{prf}
(a) We first recall that norms on a fixed finite dimensional $K$-vector space are all equivalent (see e.g.\ the argument for \cite[Theorem 3.2]{Sch70}), and one may then use the argument for \cite[Proposition 17.4(ii)]{Sch02} to finish the proof for this part.

To be more precise, let us consider $z = \sum_{j=1}^m v_j\otimes w_j\in E\otimes_K F$ with $\{v_1,...,v_m\}$ being linearly independent.
Suppose that $z = \sum_{i=1}^r x_i\otimes y_i$ as well.
The equivalent of norms as stated above gives a constant $c>0$ such that $\|\sum_{j=1}^m a_j v_j\| \geq c \cdot\max_{j=1,...,m} |a_j| \|v_j\|$, whenever $a_1,...,a_m\in K$.
As $F$ is ultrametric, it is well-known that one can find a base $\{e_1,...,e_n\}$ for the subspace of $F$ spanned by $\{w_1,...,w_m, y_1,...,y_r\}$ such that $\|\sum_{k=1}^n b_k e_k\|\geq \frac{1}{2}\max_{k=1,...,n} |b_k| \!\ \|e_k\|$ (see e.g.\ the proof of \cite[Proposition 10.4]{Sch02}).
Let $w_j = \sum_{k=1}^n a_{jk} e_k$ and $y_i = \sum_{k=1}^n b_{ik}e_k$.
Then
\begin{eqnarray*}
\sum_{i=1}^r \|x_i\|\!\ \|y_i\|
& \geq & \frac{1}{2} \sum_{i=1}^r \max_{k=1,...,n} \|x_i\|\!\ |b_{ik}| \!\ \|e_k\|
\ \geq \ \frac{1}{2} \max_{k=1,...,n} \left\|\sum_{i=1}^r b_{ik} x_i \right\| \|e_k\| \\
& \geq & \frac{c}{2} \max_{\underset{k=1,...,n}{j=1,...,m}} |a_{jk}|\!\ \|v_j\|\!\ \|e_k\|
\ \geq \ \frac{c}{2} \max_{j=1,...,m} \|v_j\| \!\ \|w_j\|
\end{eqnarray*}
Thus, if $\|z\|_\wedge = 0$, then $w_j =0$ for all $j=1,...,m$ and $z = 0$.

\smnoind
(b) $(1)\Rightarrow (2)$.
Suppose that $x\in E$ such that $\phi(x) \neq 0$.
It is clear that
\begin{equation}\label{eqt:ultra-metric-ten}
\|x\otimes 1\|_\wedge^\um = \inf\Big\{\max_{i=1,...,n} \|x_i\|: x_1,...,x_n\in E; x = \sum_{i=1}^n x_i\Big\},
\end{equation}
and we have $\|x\otimes 1\|_\wedge^\um\geq |\phi(x)| > 0$ as required.

\noindent
$(2)\Rightarrow (3)$.
As $E\hat\otimes_K^\um K \neq (0)$, there is $x\in E$ with $x\otimes 1$ being non-zero in $E\hat \otimes_K^\um K$ and hence $x\otimes 1$ is non-zero $E\hat\otimes_\BZ^\um K^0$.

\noindent
$(3)\Rightarrow (1)$.
Since $E\hat\otimes_\BZ^\um K^0$ is a non-zero ultrametric $K^0$-Banach space and $K^0$ spherically complete, the Hahn Banach theorem (see e.g.\ \cite[Proposition 9.2]{Sch02}) produces a non-zero contractive $K^0$-linear map from $E\hat\otimes_\BZ^\um K^0 = E\hat\otimes_{K^0}^\um K^0$ (see Lemma \ref{lem:proj-ten}(b)) to $K^0$ and its composition with the canonical map from $E$ to $E\hat\otimes_\BZ^\um K^0$ gives a non-zero contractive additive map.

\smnoind
(c) Let $F:= E\hat \otimes_K^\um K$ and $\iota:E = E\otimes_K K\to F$ be the canonical map.

\noindent
$\Rightarrow)$.
The hypothesis implies that $\iota$ is injective.
As $F$ is an ultrametric $K^0$-Banach space and $K^0$ is spherically complete, the Hahn Banach theorem gives a contractive additive map $\psi: F\to K^0$ with
$|\psi(\iota(x))| > 0$ as required.

\noindent
$\Leftarrow)$.
It follows from the argument of $(1)\Rightarrow (2)$ in part (b) that $\iota$ is injective, which is equivalent to $\|\cdot\|_\wedge^\um$ being a norm on $E\otimes_K K$.
\end{prf}

\medskip

By part (b) above, $E\hat\otimes_K^\um K\neq (0)$ when $E$ admits a linearly independent subset $\{e_i\}_{i\in I}$ with its linear span being norm-dense in $E$ such that there is $\kappa > 0$ with 
\begin{equation}\label{eqt:base}
\left\|\sum_{k=1}^n \lambda_k e_{i_k}\right\| \geq \kappa \cdot \max_{k=1,...,n} |\lambda_k|
\end{equation}
for any $\lambda_1,...,\lambda_n\in K$ and any distinct elements $i_1,...,i_n$ in $I$.
In particular, $E\hat\otimes_K^\um K\neq (0)$ when $E$ is finite dimensional.

\medskip

One may use projective tensor product to give a description for Banach subrings of Banach algebras.

\medskip

\begin{prop}\label{prop:reg-ring-ten-prod}
We say that a unital Banach ring $R$ is \emph{regular} if $\|m\cdot a\| = \|m\cdot 1\|\!\ \|a\|$ for any $m\in \BZ$ and $a\in R$.

\smnoind
(a) $R$ is regular if and only if there exists $K_R\in \CK^{\min}$ such that %$R\hat \otimes_\BZ K_R\neq (0)$ and
the canonical map $\Psi_R\in \CC(R;R\hat \otimes_\BZ K_R)$ is an isometry.

\smnoind
(b) Suppose that $R$ is ultrametric.
Then $R$ is regular if and only if there is $K_R\in \CK^{\um,\min}$ such that %$R\hat \otimes_\BZ^\um K_R\neq (0)$ and
the canonical map $\Phi_R\in \CC(R;R\hat \otimes_\BZ^\um K_R)$ is an isometry.
\end{prop}
\begin{prf}
The sufficiencies of both parts (a) and (b) are clear.
Assume that $R$ is regular and let $\varphi_R: \BZ\to R$ be the map $\varphi_R(m):= m\cdot 1$.

Suppose that the subring $\varphi_R(\BZ)$ is finite.
Then $\varphi_R(\BZ)$ is isomorphic to $\BZ / n\BZ$ for some $n\in \BN$.
As $k\mapsto \|k \cdot 1\|$ induces a multiplicative norm on $\BZ /n \BZ$, we know that $n\in \BP$ and $R$ is a unital $\BZ(n)$-Banach algebra.
If we take $K_R:=\BZ(n)$, then Lemma \ref{lem:univ-prop-proj-ten}(c) tells us that $\id\otimes \varphi_R$ induces a surjection in $\CC(R\hat \otimes_\BZ K_R;R)$ (respectively, $\CC(R\hat \otimes_\BZ^\um K_R;R)$) and hence the map $\Psi_R$ (respectively, $\Phi_R$) is a bijective isometry (respectively, when $R$ is ultrametric).

In the following, we consider the case when $\varphi_R(\BZ)$ is infinite, i.e.\  $\varphi_R$ is injective.
The norm on $\BZ$ induced by $\varphi_R$ produces a multiplicative norm on $\BQ$, and we set $K_R\in \CK^{\min}$ to be the completion under this norm.

\smnoind
(a) Let us first define a semi-norm $\|\cdot\|_*$ on $R\otimes_\BZ \BQ$ by
$$\|z\|_* := \inf \left\{\sum_{i=1}^l \|a_i\| |r_i|_{K_R}: l\in \BN; a_1,...,a_l\in R; r_1,...,r_l\in \BQ; z = \sum_{i=1}^l a_i\otimes r_i\right\},$$
and let $A$ be the completion of $R\otimes_\BZ \BQ/\ker \|\cdot\|_*$.
For any $b\in R$ and $\epsilon >0$, there exists $a_1,...,a_l\in R$ and $r_1,...,r_l\in \BQ$ such that $b\otimes 1 = \sum_{i=1}^l a_i\otimes r_i$ and $\sum_{i=1}^l \|a_i\| |r_i|_{K_R} \leq \|b\otimes 1\|_* + \epsilon$.
Suppose that $r_i = m_i/n_i$ with $m_i\in \BZ$ and $n_i\in \BN$ ($i=1,...,l$).
Since $R$ is regular and $\varphi_R$ is injective, $R$ is a torson-free $\BZ$-module.
Thus, $n_1\cdots n_l  b = \sum_{i=1}^l m_ia_i$.
This, together with the regularity of $R$, gives
\begin{equation}\label{eqtnorm-b-ten-1}
\|b\|
\ \leq\ \sum_{i=1}^l \|m_i\cdot 1\|\!\ \|a_i\|/\|n_i\cdot 1\|
\ =\ \sum_{i=1}^l \|a_i\|\!\  |r_i|_{K_R}
\ \leq\ \|b\otimes 1\|_* + \epsilon,
\end{equation}
which shows that $\|b\otimes 1\|_* = \|b\|$.
Thus, $b\mapsto b\otimes 1$ induces an isometry $\Psi_R\in \CC(R;A)$.
As $\|1\otimes r\|_*\leq |r|_{K_R}$ ($r\in \BQ$), the map $r\mapsto 1\otimes r$ extends to an element in $\Psi_{K_R}\in \CC(K_R;A)$.
By Lemma \ref{lem:univ-prop-proj-ten}(c), there exists $\Psi\in \CC(R\hat\otimes_\BZ K_R; A)$ such that $\Psi(b\otimes s) = \Psi_R(b)\Psi_{K_R}(s)$.
Consequently, $\|b\| = \|\Psi_R(b)\| = \|\Psi(b\otimes 1)\| \leq \|b\otimes 1\|_\wedge$ as required.

\smnoind
(b) This follows from a similar argument as part (a).
\end{prf}

\medskip

\begin{eg}\label{eg:zero-ten-prod}
Let $\BZ_1$  be the ring $\BZ$ equipped with the Euclidean norm.

\smnoind
(a) One may use Lemma \ref{lem:univ-prop-proj-ten}(c) to verify that $\BR_1\hat \otimes_\BZ \BZ_1\cong \BR_1$.
However, one has $\BR_1\hat\otimes_\BZ^\um \BZ_1 = (0)$, since $\|x\otimes 1\|_\wedge^\um$ is as given in \eqref{eqt:ultra-metric-ten}, for any $x\in \BR$.

\smnoind
(b) By Lemma \ref{lem:proj-ten}(a), we know that $\BQ_p^1\hat\otimes_\BZ \BR_1 = (0)$ (note that the image of $\BZ$ in $\BQ_p^1$ is bounded while its image in $\BR_1$ is unbounded).

\smnoind
(c) If $K$ and $L$ are (respectively, non-Archimedean) complete valuation fields with $L$ being an extension of $K$, then $K\hat \otimes_K L \cong L$ (resepectively, $K\hat \otimes^\um_K L \cong L$).

\smnoind
(d) Let $R$ be a unital ultrametric Banach ring and $K\in \CK^{\min}$ with characteristic zero.
Note that any element in $R\otimes_\BZ \BQ$ is of the form $x\otimes \frac{1}{n}$ ($x\in R;n\in \BN$).
Set
$$\left\|x\otimes \frac{1}{n}\right\| := \inf \left\{\frac{\|y\|_R}{|m|_K}: y\in R; m\in \BN \text{ with } kmx = kny \text{ for some }k\in \BZ\setminus \{0\}\right\}.$$
It is not hard to check that $\|x\otimes \frac{1}{n}\| = \inf \{\max_{i=1}^n \|x_i\|_R |r_i|_K: x_i\in R;r_i\in \BQ;x\otimes \frac{1}{n} = \sum_{i=1}^n x_i\otimes r_i\}$.
If, in addition, $R$ is regular, then $R\hat\otimes_\BZ^\um K$ is the completion of $R\otimes_\BZ \BQ$ under the above norm (see the proof of Proposition \ref{prop:reg-ring-ten-prod}), and we have $R\hat\otimes_\BZ^\um K = (0)$ when $K \neq K_R$ (note that $\|n\cdot 1\otimes \frac{1}{n}\| = \frac{|n|_{K_R}}{|n|_{K}}$ because $R$ is regular).
\end{eg}

\medskip

In the remainder of this section, we recall some facts concerning commutative rings.
Until the end of this section, $T$ is a unital commutative ring.
We set
$$\KM^0(T) := \big\{\lambda: \lambda \text{ is a non-zero multiplicative semi-norm on }T\big\}.$$
Consider $\CT_T$ to be the topology  on $\KM^0(T)$ given by pointwise convergences on elements of $T$.

\medskip

If $R$ is a commutative unital Banach ring, we denote, as in the literature,
$$\BA_R^1 := \big\{\lambda\in \KM^0(R[\bt]) :\lambda(x)\leq \|x\|, \text{ for every }x\in R\big\}$$
and call it the \emph{one-dimensional affine analytic space over $R$}.
Note that 
$$\AZ = \KM^0(\KPI).$$
Furthermore, if $R$ is ultrametric, then elements in $\BA_R^1$ are ultrametric (see Lemma \ref{lem:k-lin}(c)).
In particular, $\AZU$ consists of all ultrametric semi-norms in $\AZ$, where $\BZ_0$ is the ring $\BZ$ equipped with the trivial valuation.

\medskip

As in \cite{Berk90}, when $T$ is a Banach ring, we denote by $\KM(T)\subseteq \KM^0(T)$ the subset consisting of contractive semi-norms.
The following proposition contains some well-known facts concerning $\KM^0(T)$ (see Theorem 1.2.1 and Remark 1.2.2(i) of \cite{Berk90}).

\medskip

\begin{prop}\label{prop:mult-norm>field}
Let $T$ be a commutative unital ring and $\lambda\in \KM^0(T)$.

\smnoind
(a) There is a complete valuation field $\CH(\lambda)$ and a unital ring homomorphism $\varphi_\lambda:T\to \CH(\lambda)$ such that $\CH(\lambda)$ is generated by $\varphi_\lambda(T)$ as a complete valuation field and $\lambda(r) = |\varphi_\lambda(r)|$ ($r\in T$).
Moreover, $(\CH(\lambda),\varphi_\lambda)$ is unique up to isometric isomorphism.

\smnoind
(b) If, in addition,  $T$ is a Banach ring, then $\KM(T)$ is a non-empty compact Hausdorff space.
\end{prop}
%\begin{comment}
[This proof will not appear in the published version.]

\begin{prf}
(a) As $\ker\lambda$ is a prime ideal of $T$, the quotient $T/\ker \lambda$ is an integral domain.
We denote by $F$ its field of quotients. 
Since $\lambda$ induces a multiplicative norm on $T/\ker \lambda$, it induces a valuation on $F$.
We denote by $\CH(\lambda)$ its completion and by $\varphi_\lambda$ the canonical embedding. 
The uniqueness of $(\CH(\lambda),\varphi_\lambda)$ is clear. 

\smnoind
(b) This is precisely \cite[Theorem 1.2.1]{Berk90}. 
\end{prf}
%\end{comment}

\medskip

Let us also present the following example.
Notice that parts (a) and (b) are adapted from \cite[Remark 2.2.2]{Berk90}.

\medskip

\begin{eg}\label{eg:aff-neigh}
Consider $M\in \BN$ and pick any $\bp_1,...,\bp_m,\bq_1,...,\bq_m\in \KPI$ as well as $\upsilon_1,...,\upsilon_m,\omega_1,...,\omega_m > 0$.
If $R$ is a unital Banach ring, we denote by $R\la \upsilon_1^{-1}\bt_1,\omega_1\bt_2,...,\upsilon_m^{-1}\bt_{2m-1},
\omega_m\bt_{2m}\ra$ the Banach ring  of formal power series satisfying certain norm summability properties as in \cite[Example 1.1.1(v)]{Berk90}.
In the case when $R$ is ultrametric, we denote by $R\{ \upsilon_1^{-1}\bt_1,\omega_1\bt_2,...,\upsilon_m^{-1}\bt_{2m-1},
\omega_m\bt_{2m}\}$ the Banach ring  of formal power series satisfying certain norm convergence properties as in \cite[\S 2.1]{Berk90}.

\smnoind
(a) Let $\BZ_1$ be as in Example \ref{eg:zero-ten-prod}.
Since $M\geq 1$, the Banach ring $\BZ_1\la M^{-1} \bt\ra$ coincides with $\KPI$ as unital rings and
\begin{equation}\label{eqt:spec-Z-t}
\KM(\BZ_1\la M^{-1}\bt\ra) = \{\gamma\in \AZ: \gamma(\bt)\leq M\}.
\end{equation}
Let us denote by $\BZ_1\la M^{-1} \bt\ra_{\upsilon^{-1} \bp, \omega\bq^{-1}}$ the quotient of the unital Banach ring \linebreak
$\BZ_1\la M^{-1} \bt \ra \la \upsilon_1^{-1}\bt_1,\omega_1\bt_2,...,\upsilon_m^{-1}\bt_{2m-1},
\omega_m\bt_{2m}\ra$ under the closed ideal generated by $\{\bt_{2i-1}-\bp_i\}_{i\in \{1,...,m\}}\cup \{\bq_j\bt_{2j}-1\}_{j\in \{1,...,m\}}$.
The canonical map $\tau: \BZ_1\la M^{-1} \bt\ra\to \BZ_1\la M^{-1} \bt\ra_{\upsilon^{-1} \bp, \omega\bq^{-1}}$ induces a homeomorphism %$\ti\tau$
from $\KM(\BZ_1\la M^{-1}\bt\ra_{\upsilon^{-1} \bp, \omega\bq^{-1}})$ onto
$$\big\{\lambda\in \AZ: \lambda(\bt)\leq M; \lambda(\bp_i) \leq   \upsilon_i, \lambda(\bq_i) \geq \omega_i \text{ for } 1\leq i \leq m\big\}.$$

\smnoind
(b) Let $\BZ_0$ be as in the second paragraph after Example \ref{eg:zero-ten-prod}.
As in part (a), $\BZ_0\{M^{-1}\bt\}$ coincides with $\KPI$ as rings and
$$\KM(\BZ_0\{ M^{-1}\bt\}) = \{\gamma\in \AZU: \gamma(\bt)\leq M\}.$$
Suppose that $\BZ_0\{M^{-1}\bt\}_{\upsilon^{-1} \bp, \omega\bq^{-1}}$ is the quotient of
the ultrametric Banach ring $\BZ_0\{ M^{-1} \bt\}\{\upsilon_1^{-1}\bt_1,\omega_1\bt_2,...,\upsilon_m^{-1}\bt_{2m-1},
\omega_m\bt_{2m}\}$ under the closed ideal generated by $\{\bt_{2i-1}-\bp_i\}_{i\in \{1,...,m\}}\cup \{\bq_j\bt_{2j}-1\}_{j\in \{1,...,n\}}$.
%with the quotient map $Q^\um_{\upsilon^{-1} \bp, \omega\bq^{-1}}$.
The canonical map
$\tau^\um: \BZ_0\{ M^{-1}\bt\}\to \BZ_0\{ M^{-1}\bt\}_{\upsilon^{-1} \bp, \omega\bq^{-1}}$ induces a homeomorphism from $\KM(\BZ_0\{ M^{-1}\bt\}_{\upsilon^{-1} \bp, \omega\bq^{-1}})$ onto
$$\big\{\lambda\in \AZU: \lambda(\bt)\leq M; \lambda(\bp_i) \leq   \upsilon_i, \lambda(\bq_i) \geq \omega_i \text{ for } 1\leq i \leq m\big\}.$$

\smnoind
(c) Let $\lambda\in \KM(\BZ_1\la M^{-1}\bt\ra)$.
There exist sequence $\{\bp_i\}_{i\in \BN}$ and $\{\bq_i\}_{i\in \BN}$ in $\KPI$ as well as sequences $\{\upsilon_i\}_{i\in \BN}$ and $\{\omega_i\}_{i\in \BN}$ in $(0,\infty)$ satisfying the following conditions:
\begin{enumerate}[1).]
\item if $\lambda(\bp_i) >0$, then one has $\bq_i = \bp_i$ and $\omega_i < \lambda(\bp_i) < \upsilon_i$;
\item if $\lambda(\bp_i) = 0$, then one has $\bq_i = 1$, $\omega_i = 1/2$  and $\lambda(\bp_i) < \upsilon_i$;
\item for any $\bp\in \KPI$ and $\epsilon >0$, there is $i \in \BN$ with $\bp_i = \bp$, $\upsilon_i \leq \lambda(\bp)+\epsilon$ and $\omega_i \geq \lambda(\bp) -\epsilon$.
\end{enumerate}
Let $T:= \BZ_1 \la M^{-1} \bt \ra \la \upsilon_1^{-1}\bt_1,\omega_1\bt_2,\upsilon_2^{-1}\bt_3,\omega_2\bt_4, ...\ra$ (which is defined in a similar way as the Banach ring in \cite[Example 1.1.1(v)]{Berk90}), and $I$ be the  closed ideal generated by $\{\bt_{2i-1}-\bp_i\}_{i\in \BN}\cup \{\bq_j\bt_{2j}-1\}_{j\in \BN}$.
Then $T / I \cong \CH(\lambda)$.
A similar fact holds for $\lambda\in \KM(\BZ\{ M^{-1}\bt\})$.

\smnoind
(d) If $\Gamma$ is a (discrete) abelian group, $\ell^1(\Gamma;\BC)$ is the unital Banach ring of absolutely summable complex functions on $\Gamma$ and $\hat \Gamma$ is the dual (compact) group of $\Gamma$, then $\KM(\ell^1(\Gamma;\BC)) = \hat \Gamma/\sim_{\rm cong}$, where $f\sim_{\rm cong} g$ if and only if $f(t) \in \{g(t),\overline{g(t)}\}$ ($t\in \Gamma$).
\end{eg}

\medskip

We end this preliminary section with the following well-known proposition.
We give a brief account here for the sake of completeness.

\medskip

\begin{prop}\label{prop:top-AZ}
	$\AZ$ is a second countable locally compact Hausdorff space under the topology $\CT_\KPI$ and $\AZU$ is a closed subset of $\AZ$.
	Moreover, for any $M \in \BN$, the subset
	$\{\lambda\in \AZ: \lambda(\bt) \leq M\}$
	is compact. 	
\end{prop}
\begin{prf}
	The compactness of $\{\lambda\in \AZ: \lambda(\bt) \leq M\}$ and the local compactness of $\AZ$ follows from Relation \eqref{eqt:spec-Z-t}, Proposition \ref{prop:mult-norm>field}(b) as well as
	$$\AZ = \bigcup_{N\in \BN} \{\lambda\in \AZ: \lambda(\bt) < N\}.$$
	The second countability follows from the fact that for any $\lambda\in \AZ$, the collection of subsets of the form
\begin{equation}\label{eqt:neigh-TBK}
U_{\upsilon^{-1}\bp,\omega \bq^{-1}}:= \big\{\gamma\in \AZ: \gamma(\bp_i) <   \upsilon_i, \gamma(\bq_i) > \omega_i \text{ for } 1\leq i \leq m\}
\end{equation}
%changed after submission
	(where $\bp_1,...,\bp_m, \bq_1, ...,\bq_m\in \KPI$ and $\upsilon_1,...,\upsilon_m , \omega_1,...,\omega_m \in \BQ_+\setminus \{0\}$) is a base of open neighborhoods of for the topology $\CT_\KPI$ (see e.g.\ \cite[Remark 2.2.2]{Berk90}).
	The closedness of $\AZU$ follows from Lemma \ref{lem:k-lin}(c).
	
\end{prf}

\medskip

\section{The universal space of scalars}

From now on, $\alpha$ is an infinite cardinal and $\CK$ (respectively, $\CK^\um$) is the set of all complete valuation fields (respectively, non-Archimedean complete valuation fields) with cardinality not exceeding $2^\alpha$.

\medskip

If $K_1,K_2\in \CK$ and $s_i\in K_i$ ($i=1,2$), we define  $(K_1,s_1)\sim (K_2,s_2)$ whenever there exist $K\in \CK$, $s\in K$ as well as $\varphi_i\in \CC(K;K_i)$ with $\varphi_i(s) = s_i$ ($i=1,2$).
As $\varphi_1$ and $\varphi_2$ are isometric (by Lemma \ref{lem:k-lin}(a)), $\sim$ is an equivalence relation on the disjoint union of all elements in $\CK$.
We denote by $(K_1,s_1)_\sim$ the equivalence class containing $(K_1,s_1)$ and by $\ti \BK$ the set of all equivalence classes.
Observe that if $K_1$ is non-Archimedean, then so is $K_2$ and we denote by $\ti \BK^{\um}$ the subset of $\ti \BK$ consisting of equivalence class of elements in non-Archimedean fields.
Moreover, we set $\mu_{L}$ to be the map from $L\in \CK$ to $\ti \BK$ that sends $r\in L$ to $(L,r)_\sim$.
Clearly, if $K\in \CK$ is an extension of $L$, then $\mu_{K}(r) = \mu_{L}(r)$.

\medskip

Consider $L^{r}$ to be the closed subfield of $L$ generated by $r\in L$.
In particular, $L^0 = L^1$ is the smallest closed subfield of $L$.
Observe that $(L,r)\sim (L^{r},r)$.
Moreover, $(K_1,s_1)\sim (K_2,s_2)$ if and only if there is $\varphi\in \CC(K_1^{s_1};K_2^{s_2})$ with $\varphi(s_1) = s_2$.

\medskip

Suppose that $L\in \CK$ is singly generated and $G_L$ is the set of generators of $L$.
For any $r,s\in G_L$, one has $(L,r)\sim (L,s)$ if and only if there is $\varphi\in \Aut(L)$  such that $s = \varphi(r)$.
Therefore, $\mu_L$ induces an injection from $G_L/\Aut(L)$ to $\ti \BK$.
Moreover, $\ti \BK$ is the union of images of $G_L/\Aut(L)$, when $L$ runs through all singly generated fields.

\medskip

\begin{rem}\label{rem:K-ti-s}
(a) Suppose that $(K_1,s_1)\sim (K_2,s_2)$.
If $K,s,\varphi_1$ and $\varphi_2$  are as in the above, then for any $\bp\in \KPI$, one has $\varphi_i(\bp(s)) = \bp(s_i)$ ($i=1,2$) and $|\bp(s_1)|_{K_1} = |\bp(s)|_K = |\bp(s_2)|_{K_2}$.
Therefore, any $(L,r)_\sim\in \ti \BK$ induces an element $\lambda_{(L,r)_\sim}$ in $\AZ$ with
$$\lambda_{(L,r)_\sim}(\bp) := |\bp(r)|_L \qquad (\bp\in \KPI).$$
Notice that if $L\in\CK^\um$, then $\lambda_{(L,r)_\sim}(\bp)\in \AZU$.

\smnoind
(b) For any $(K,s)_\sim\in \ti \BK$, it is easy to see that the field $\CH(\lambda_{(K,s)_\sim})$ as in Proposition \ref{prop:mult-norm>field}(a) is the unique complete valuation field with a generator $r$ such that  $\big(\CH(\lambda_{\ti s}),r\big)\sim (K,s)$  (actually, $r=\varphi_{\lambda_{(K,s)_\sim}}(\bt)$).
Thus, one may find a surjective isometry $\Phi\in \CC\big(\CH(\lambda_{\ti s});K^s\big)$ with $\Phi\big(\varphi_{\lambda_{(K,s)_\sim}}(\bt)\big) = s$.

\smnoind
(c) For the sake of simplicity, we will sometimes denote $(K,s)_\sim$ by $\ti s$.
\end{rem}

\medskip

The following proposition could be known but since we do not find it in the literature, we present a proof here. 

\medskip

\begin{prop}\label{prop:top-K-alpha}
(a) The map $\Theta:\ti s \mapsto \lambda_{\ti s}$ is a bijection from $\ti \BK$ onto $\AZ$ such that $\Theta \big(\ti \BK^{\um }\big) = \AZU$.

\smnoind
(b) $\Theta \circ \mu_K$ is continuous for any $K\in \CK$.

\smnoind
(c) If $K\in \CK^{\min}$, then $\Theta \circ \mu_K:K\to \Theta (\mu_K(K))$ is a homeomorphism.

\smnoind
(d) $\mu_K(K)$ is a closed subset of $\ti \BK$ when $K\in \CK^{\min}\setminus \{\BQ_0\}$.
\end{prop}
\begin{prf}
(a) Suppose that $(K_1,s_1)_\sim,(K_2,s_2)_\sim\in \ti \BK$ with $\lambda_{(K_1,s_1)_\sim} =
\lambda_{(K_2,s_2)_\sim}$.
For each $\bp\in \KPI$, one has $|\bp(s_1)|_{K_1} = |\bp(s_2)|_{K_2}$.
Hence, if $\bp,\bq,\bp ',\bq '\in \KPI$ satisfying $\bp(s_1)\neq 0$, $\bq(s_1)\neq 0$ and $\bp '(s_1)\bp(s_1)^{-1} = \bq '(s_1)\bq(s_1)^{-1}$, then $\bp(s_2)\neq 0$, $\bq(s_2)\neq 0$ and $(\bp ' \bq - \bq '\bp)(s_2) = 0$.
This gives a well-defined isometric ring homomorphism from the (not necessarily closed) subfield of $K_1$ generated by $s_1$ to $K_2^{s_2}$  that sends $\bp ' (s_1) \bp(s_1) ^{-1}$ to $\bp '(s_2)\bp(s_2)^{-1}$.
It extends to an
isometry $\varphi\in \CC(K_1^{s_1};K_2^{s_2})$.
In particular,  $(K_1,s_1)\sim (K_2,s_2)$.
Hence, $\Theta $ is injective.

To show the surjectivity, we consider $\lambda\in \AZ$ and let  $(\CH(\lambda),\varphi_\lambda)$ be as in Proposition \ref{prop:mult-norm>field}(a).
It is not hard to see that $\lambda = \Theta \big((\CH(\lambda), \varphi_\lambda(\bt))_\sim\big)$.
Finally, it is easy to verify that $\Theta \big(\ti \BK^{\um }\big) = \AZU$.

\smnoind
(b) This part is clear.

\smnoind
(c) To simplify the notation, we will ignore $\Theta $ and identify $\ti \BK$ with $\AZ$ directly.
If $(K,r)\sim (K,s)$, then the minimality of $K$ implies that $K^r = K = K^s$, and Lemma \ref{lem:k-lin}(b) gives $r=s$.
This shows that $\mu_K$ is injective.
Next, we show that if $\{s_i\}_{i\in \BN}$ is a sequence in $K$ and $s\in K$ with $\mu_K(s_i) \to \mu_K(s)$, then $|s_i - s|_K \to 0$.

Let us first consider the case when $K\in \CK^{\min}_\BR$ or $K\in \CK^{\min}_{\BQ_p}$ for some $p\in \BP$.
Then any closed and bounded subset of $K$ is compact.
For every $\delta >0$, we set
$$C_{\delta}:=\{r\in K: |r|_K \leq \delta + 1\}.$$
Since $|s_i|_K \to |s|_K$ (by considering $\bt\in \KPI$), we may assume that $s_i\in C_{|s|_K}$ for all $i\in \BN$.
The compactness of $C_{|s|_K}$ tells us that $\mu_K|_{C_{|s|_K}}$ is a homeomorphism, and hence $s_i\to s$.

Secondly, if $K=\BZ(p)$ for some $p\in \BP$, then $K$ is compact and the conclusion follows.

Finally, we consider the case when $K=\BQ_0$.
For any $r\in \BQ$, one has
\begin{equation}\label{eqt:val-on-BQ}
\mu_K(r)(\bp) = \begin{cases}
0\ & \text{if }r \text{ is a zero of }\bp\\
1  & \text{otherwise}
\end{cases}
\qquad (\bp\in \KPI).
\end{equation}
Suppose that $s = m/n$, where $m\in \BZ$ and $n\in \BN$.
Set $\bp_0:=n\bt - m\in \KPI$.
The fact that $\mu_K(s_i)(\bp_0) \to \mu_K(s)(\bp_0)$ implies that $\bp_0(s_i) = 0$ (or equivalently, $s_i = m/n$), when $i$ is large enough.

\smnoind
(d) Let $\{s_i\}_{i\in \BN}$ be a sequence in $K$.
Suppose that $L\in \CK$ and $r\in L$ satisfying $\mu_K(s_i) \to \mu_{L}(r)$.
As in the argument of part (d), we may assume that $s_i\in C_{|r|_L}$ for all $i\in \BN$.
The compactness of $C_{|r|_L}$ (as $K\neq \BQ_0$) produces a subsequence $\{s_{i_j}\}_{j\in \BN}$ that converges to some $s\in K$.
Thus, $(L,r)_\sim  = (K,s)_\sim \in \mu_K(K)$ by the continuity of $\mu_K$.
\end{prf}

\medskip

The above tells us that $\ti \BK$ and $\ti \BK^{\um }$ are independent of the infinite cardinal $\alpha$.
We will identify $\ti \BK$ (respectively, $\ti \BK^\um$) with the topological space $\AZ$ (respectively, $\AZU$) directly through Proposition \ref{prop:top-K-alpha}(a) (i.e.\  ignore $\Theta $).
For any $\lambda\in \AZ$, we set $|\lambda|:= \lambda(\bt)$ as well as
\begin{equation}\label{eqt:TBK_M}
(\AZ)_{M} := \{\lambda\in \TBK: |\lambda| \leq M\}
\quad \text{and} \quad
(\AZU)_M := \{\lambda\in \AZU: |\lambda| \leq M\}
\quad (M\geq 0).
\end{equation}
As stated in Proposition \ref{prop:top-AZ}, these two are compact sets.
Moreover, if $\lambda = \lambda_{(K,s)_\sim}$, then obviously, $|\lambda| = |s|_K$.

\medskip

\begin{rem}\label{rem:neigh-TBK-M}
(a) The compact set $(\AZ)_1$ contains the equivalence classes of all elements in all fields equipped with the trivial valuations.

\smnoind
(b) By the argument of Proposition \ref{prop:top-AZ} and Example \ref{eg:aff-neigh}(a), for every $M\in \BN$ and $\lambda\in \TBK$ with $\lambda(\bt) < M$, the collection of subsets of the form $\KM(\BZ_1\la M^{-1}\bt\ra_{\upsilon^{-1} \bp, \omega\bq^{-1}})$ containing $\lambda$ is a compact neighborhood base for $\lambda$.
\end{rem}

\medskip

Since $\BQ_0$ is an infinite discrete subset of the compact subset $(\AZ)_1$, $\BQ_0$ is not closed in $\TBK$, i.e.\ Proposition \ref{prop:top-K-alpha}(d) fails in the case when $K=\BQ_0$.
We will see in part (a) of the following an unusual cluster point of  $\BQ_0$ in $\TBK$.

\medskip

\begin{eg}\label{eg:non-inj}
(a) Let $\BZ(\bt)$ be the field of quotients of the integral domain $\KPI$.
We equip $\BZ(\bt)$ with the trivial valuation.
Suppose that $\{r_i\}_{i\in \BN}$ is a sequence in $\BQ_0$ such that $r_i\neq r_j$ when $i\neq j$.
Since $|\bp(r_i)|_{\BQ_0} \to 1$ for any $\bp\in \KPI\setminus \{0\}$, we know that $(\BQ_0,r_i)_\sim \to (\BZ(\bt), \bt)_\sim$ in $\TBK$.

\smnoind
(b) Consider $\BC_1$ to be the field $\BC$ equipped with the Euclidean valuation $|\cdot|_1$.
One may identify $\mu_{\BC_1}(\BC_1)$ with the upper half plane $\mathbb{H}$, under which $\mu_{\BC_1}$ is the map $Q_{\mathbb{H}}:\BC_1\to \mathbb{H}$ given by $Q_{\mathbb{H}}(s):= \begin{cases}
s \ & \text{if }s\in \mathbb{H}\\
\bar s & \text{otherwise.}
\end{cases}$

In fact, since the set $G_\BC$ of generators of $\BC_1$ equals $\BC\setminus \BR$ and $\CC({\BC_1};{\BC_1}) = \{\id; \mathfrak{c}\}$, where $\mathfrak{c}$ is the complex conjugation, we see that $\mu_{\BC_1}(r) = \mu_{\BC_1}(s)$ if and only if $r\in \{s,\bar s\}$ (see the discussion preceding Remark \ref{rem:K-ti-s}).
Thus, $\mu_{\BC_1}$ induces a continuous bijection $\check \mu_{\BC_1}: \mathbb{H} \to
\mu_{\BC_1}({\BC_1})$ with $\check \mu_{\BC_1}(Q_{\mathbb{H}}(s)) = \mu_{\BC_1}(s)$.

Now, suppose that $\{s_i\}_{i\in \BN}$ is a sequence in $\mathbb{H}$ and $s\in \mathbb{H}$ such that $\check \mu_{\BC_1}(s_i) \to \check \mu_{\BC_1}(s)$.
One may assume that every $s_i$ ($i\in \BN$) lies in the set $C_{|s|_1}:=\{r\in \mathbb{H}: |r|_1 \leq |s|_1+1\}$.
As $C_{|s|_1}$ is compact and $\check \mu_{\BC_1}$ is a continuous injection, we know that $s_i\to s$.

\smnoind
(c) For any $s\in \BC$, we set $\KPI_s := \{\bp\in\KPI: \bp(s) = 0\}$. 
If $r,s\in \BC$, we define $r\sim_0 s$ if $\KPI_r = \KPI_s$. 
It is easy to see that $\mu_{\BC_0}$ induces an injection from $\BC/\sim_0$ to $\TBK$. 
If $\{s_i\}_{i\in \KI}$ is a net in $\BC$ and $s\in \BC$, then $\mu_{\BC_0}(s_i)\to \mu_{\BC_0}(s)$ if and only if $\lim_i \KPI_{s_j} = \KPI_{s}$, namely,
$$\bigcap_{i\in \KI}\bigcup_{j\geq i} \KPI_{s_j} \subseteq \KPI_s \subseteq \bigcup_{i\in \KI}\bigcap_{j\geq i} \KPI_{s_j}.$$ 
\end{eg}

\medskip

In the remainder of this section, we will give a complete description of the topological subspace:
$$\TBKM:= \bigcup_{K\in \CK^{\min}} \mu_K(K).$$

\medskip

\begin{lem}\label{lem:top-K-min}
(a) $\mu_{\BQ_0}(\BQ_0)$ is closed in $\TBKM$.

\smnoind
(b) For any $K, L\in \CK^{\min}$ with $K\neq L$, one has $\mu_K(K)\cap \mu_L(L) = \emptyset$.

\smnoind
(c) $\bigcup_{\upsilon\in(0,1]} \mu_{\BR_\upsilon}(\BR_\upsilon)\cong (0,1]\times \BR_1$ as topological spaces and is an open subset of $\TBKM$ with its closure in $\TBKM$ being $\mu_{\BQ_0}({\BQ_0})\cup \bigcup_{\upsilon\in(0,1]} \mu_{\BR_\upsilon}(\BR_\upsilon)$.

\smnoind
(d) For each $q\in \BP$, the subset $\bigcup_{\omega\in (0,\infty)} \mu_{\BQ_q^\omega}(\BQ_q^\omega)$ is open in $\TBKM$ and is homeomorphic to $(0,\infty)\times \BQ_q^1$.

\smnoind
(e) If $q\in \BP$, the subset $\mu_{\BZ(q)}(\BZ(q))$ is contained in the closure of $\bigcup_{\omega\in (0,\infty)} \mu_{\BQ_q^\omega}(\BQ_q^\omega)$.

\smnoind
(f) $\mu_{\BQ_0}({\BQ_0})$ is contained in the closure of $\bigcup_{p\in \BP} \mu_{\BQ_p^1}(\BQ_p^1)$.

\smnoind
(g) The closure of $\bigcup_{p\in \BP} \bigcup_{\omega\in (0,\infty)} \mu_{\BQ_p^\omega}(\BQ_p^\omega)$ in $\TBKM$ equals
$$\BA_{\BZ_0}^{1,\min}:= \mu_{\BQ_0}({\BQ_0}) \cup \bigcup_{p\in \BP} \mu_{\BZ(p)}(\BZ(p)) \cup \bigcup_{p\in \BP}\bigcup_{\omega\in (0,\infty)} \mu_{\BQ_p^\omega}(\BQ_p^\omega).$$

\smnoind
(h) For each $q\in \BP$, the set $\mu_{\BQ_0}(\BQ_0)$ is contained in the closure of $\bigcup_{\omega\in (0,\infty)} \mu_{\BQ_q^\omega}(\BQ_q^\omega)$.
\end{lem}
\begin{prf}
(a) Suppose that $\{r_i\}_{i\in \BN}$ is a sequence in ${\BQ_0}$ and there exist $K\in \CK^{\min}$ as well as $s\in K$ with $\mu_{\BQ_0}(r_i) \to \mu_{K}(s)$.
For any $n\in \BZ\subseteq \KPI$, we know that $\mu_{\BQ_0}(r_i)(n) = |n|_{\BQ_0}$, and hence $|n|_{K} = |n|_{\BQ_0}$.
Now, we conclude from the minimality of $K$ that $K = {\BQ_0}$.

\smnoind
(b) Suppose that there exist $r\in K$ and $s\in L$ with $(K,r)\sim (L,s)$.
The minimality of $K$ and $L$ gives a bijection $\varphi\in \CC(K;L)$ such that $\varphi(r) =s$, and we have a contradiction that $K=L$.

\smnoind
(c) Let us first show that the complement of $\bigcup_{K\in \CK^{\min}_\BR} \mu_K(K)$ is closed in $\TBKM$.
Suppose that $\{L_i\}_{i\in \BN}$ is a sequence in $\CK^{\min} \setminus  \CK^{\min}_\BR$, $s_i\in L_i$ ($i\in \BN$), $K\in \CK^{\min}$ and $r\in K$ such that $\mu_{L_i}(s_i)\to \mu_K(r)$.
Then, for any $n\in \BZ$, one has $|n|_K = \mu_K(r)(n) = \lim_i |n|_{L_i}\leq 1$, which implies that $K$ is non-Archimedean.
Consequently, $K\notin \CK^{\min}_\BR$.

Next, we prove that $\mu_{\BQ_0}({\BQ_0})\cup \bigcup_{K\in \CK^{\min}_\BR} \mu_K(K)$ is closed  in $\TBKM$.
Suppose that $\{K_i\}_{i\in \BN}$ is a sequence in $\CK^{\min}_\BR \cup \{{\BQ_0}\}$, $r_i\in K_i$ ($i\in \BN$), $L\in \CK^{\min}$ and $s\in L$ such that $\mu_{K_i}(r_i)\to \mu_L(s)$.
Again, we have $|n|_{K_i} \to |n|_L$ ($n\in \BZ$).
Since $K_i\in \CK^{\min}_\BR\cup \{{\BQ_0}\}$, there is $\upsilon_i\in [0,1]$ such that $|r|_{K_i} = |r|_\E^{\upsilon_i}$ ($r\in K_i$) (and we identify $K_i\subseteq \BR$ directly).
Furthermore, as $\{\upsilon_i\}_{i\in \BN}$ is a sequence in $[0,1]$, it has a subsequence $\{\upsilon_{i_j}\}_{j\in \BN}$ that converges to some $\upsilon\in [0,1]$.
Thus, $|n|_{K_{i_j}} \to |n|_\E^\upsilon$, which implies that $|n|_L = |n|_\E^\upsilon$ ($n\in \BZ$).
This, and the minimality of $L$, tells us that $L\in \CK^{\min}_\BR\cup \{{\BQ_0}\}$.

We now verify that $\bigcup_{K\in \CK^{\min}_\BR} \mu_K(K)$ is dense in $\mu_{\BQ_0}({\BQ_0})\cup \bigcup_{K\in \CK^{\min}_\BR} \mu_K(K)$.
Suppose that $r\in \BQ$.
For each $n\in \BN$, take $r_n = r\in \BR$.
Then $\mu_{\BR^{1/n}}(r_n)(\bp) = |\bp(r)|_\E^{1/n} \to |\bp(r)|_\E^0 = \mu_{\BQ_0}(r)(\bp)$ ($\bp\in \KPI$), i.e.\ $\mu_{\BR^{1/n}}(r_n) \to \mu_{\BQ_0}(r)$ as required.

Finally, we will verify that the canonical map $\Phi: (0,1]\times \BR \to \bigcup_{\upsilon\in(0,1]} \mu_{\BR_\upsilon}(\BR_\upsilon)$ that sends $(\upsilon, t)$ to $\mu_{\BR_\upsilon}(t)\in \mu_{\BR_\upsilon}(\BR_\upsilon)$ is a homeomorphism.

In fact, $\Phi$ is bijective because of part (b) as well as Proposition \ref{prop:top-K-alpha}(c).
Let $(\upsilon, t) \in (0,1]\times \BR_1$ and $\{(\upsilon_i, t_i)\}_{i\in \BN}$ be a sequence in $(0,1]\times \BR_1$.
If $\upsilon_i\to \upsilon$ and $t_i\to t$, then $\Phi(\upsilon_i,t_i)(\bp) = |\bp(t_i)|_\E^{\upsilon_i} \to |\bp(t)|_\E^{\upsilon} = \Phi(\upsilon,t)(\bp)$ ($\bp\in \KPI$).
Conversely, suppose that $\Phi(\upsilon_i,t_i)\to \Phi(\upsilon, t)$.
Then $2^{\upsilon_i} = |2|_\E^{\upsilon_i} \to 2^\upsilon$ (as $2\in \BZ\subseteq \KPI$) and we have $\upsilon_i\to \upsilon$.
Moreover, we have $|t_i|_\E^{\upsilon_i}\to |t|_\E^\upsilon$ (by considering $\bt\in \KPI$), which implies that $\upsilon_i \ln (|t_i|_\E) \to \upsilon \ln (|t|_\E)$ (we take $\ln 0 = -\infty$ as usual).
Thus, $|t_i|_\E \to |t|_\E$.
Therefore, one can find $N\in \BN$ such that $\frac{1}{N} < \upsilon$ and
$(\upsilon_i, t_i)\in [\upsilon - \frac{1}{N},1]\times \{r\in \BR: |r|_\E\leq N\}$.
As the restriction of $\Phi$ on the compact set $[\upsilon - \frac{1}{N},1]\times \{r\in \BR: |r|_\E\leq N\}$ is a homeomorphism, we know that $(\upsilon_i,t_i)\to (\upsilon, t)$ as required.

\smnoind
(d) We first show that $\bigcup_{p\in \BP} \bigcup_{K\in \CK^{\min}_{\BQ_p}} \mu_K(K)$ is open in $\TBKM$.
Suppose on the contrary that this set is not open.
As $\bigcup_{p\in \BP} \mu_{\BZ(p)}(\BZ(p)) \cup \bigcup_{p\in \BP}\bigcup_{K\in \CK^{\min}_{\BQ_p}} \mu_K(K)$ is open in $\TBKM$ (by part (c)), there exist $q\in \BP$, $\omega\in (0,\infty)$,
$s\in \BQ_{q}^{\omega}$, a sequence $\{p_i\}_{i\in \BN}$ in $\BP$ as well as $r_i\in \BZ(p_i)$ ($i\in \BN$) such that $\mu_{\BZ(p_i)}(r_i) \to \mu_{\BQ_{q}^{\omega}}(s)$.
Observe that if $n_i\in \BZ$ with $r_i$ being the image of $n_i$, then
\begin{equation}\label{eqt:Z-qZ}
\mu_{\BZ(p_i)}(r_i)(\bp) = \begin{cases}
0 \ & \text{if }p_i \text{ divides }\bp(n_i)\\
1   & \text{otherwise},
\end{cases}
\qquad (\bp\in \KPI).
\end{equation}
However, this contradicts the fact that $\{|\bp(s)|_{q}^{\omega}:\bp\in \KPI\}\nsubseteq \{0,1\}$ (actually, $\{|n|_{q}^{\omega}:n\in \BZ\}\nsubseteq \{0,1\}$).

In order to verify the openness of $\bigcup_{\omega\in (0,\infty)} \mu_{\BQ_q^\omega}(\BQ_q^\omega)$, let us suppose that $\omega\in (0,\infty)$ and $s\in \BQ_q^\omega$ such that there exist a sequence $\{(q_i, \omega_i)\}_{i\in \BN}$ in $\BP\times (0,\infty)$ as well as $s_i\in \BQ_{q_i}^{\omega_i}$ ($i\in \BN$) with $\mu_{\BQ_{q_i}^{\omega_i}}(s_i) \to \mu_{\BQ_{q}^{\omega}}(s)$.
As $q\in \BZ\subseteq \KPI$, one has $|q|_{q_i}^{\omega_i} \to |q|_q^\omega = q^{-\omega}$.
However, since $|q|_{q_i}^{\omega_i} = 1$ when $q_i\neq q$, we know that $q_i = q$ when $i$ is large enough.
Consequently, $\bigcup_{\omega\in (0,\infty)} \mu_{\BQ_q^\omega}(\BQ_q^\omega)$ is open in $\TBKM$.

To show the second claim, let us consider the canonical map $\Psi_q: (0,\infty)\times \BQ_q^1\to \bigcup_{\omega\in (0,\infty)} \mu_{\BQ_q^\omega}(\BQ_q^\omega)$ sending $(\omega, s)$ to $\mu_{\BQ_q^\omega}(s)\in \mu_{\BQ_q^\omega}(\BQ_q^\omega)$.
As in part (c), the map $\Psi_q$ is continuous and bijective.
Suppose that $\{(\omega_i,s_i)\}_{i\in \BN}$ is a sequence in $(0,\infty)\times \BQ_q$ such that $\Psi_q(\omega_i,s_i)\to \Psi_q(\omega, s)$.
From $q^{-\omega_i} = |q|_q^{\omega_i} \to q^{-\omega}$, we know that $\omega_i\to \omega$.
This, together with $|s_i|_q^{\omega_i}\to |s|_q^\omega$, implies that $|s_i|_q \to |s|_q$.
Therefore, the same compactness argument as in part (c) tells us that $(\omega_i,s_i)\to (\omega, s)$.

\smnoind
(e) For any $s\in \BZ(q)$, there is $n\in \BZ$ with $s$ being the image $\bar n$ of $n$.
We claim that $\mu_{\BQ_q^{\omega_i}}(n) \to \mu_{\BZ(q)}(s)$ if $\omega_i\to \infty$.
In fact, for any $\bp\in \KPI$, one knows that $|\bp(n)|_q^{\omega_i} \to 0$ when $q$ divides $\bp(n)$ (in this case, $|\bp(n)|_q < 1$) and $|\bp(n)|_q^{\omega_i} =1$ ($i\in \BN$) when $q$ does not divides $\bp(n)$.
The conclusion now follows from \eqref{eqt:Z-qZ}.

\smnoind
(f) Consider $\{p_i\}_{i\in \BN}$ to be a sequence in $\BP$ with $p_i\to \infty$.
For any $r\in \BQ$, we set $r_i := r\in \BQ_{p_i}$ ($i\in \BN$).
Clearly, when $\bp\in \KPI$, one has $|\bp(r_i)|_{p_i} = 0$ if $\bp(r) = 0$ and $|\bp(r_i)|_{p_i} \to 1$ if $\bp(r)\neq 0$ (as the numerator and denominator of $\bp(r)$ have only finite numbers of prime factors).
By \eqref{eqt:val-on-BQ}, one has $\mu_{\BQ_{p_i}^1}(r_i) \to \mu_{\BQ_0}(r)$ as required.

\smnoind
(g) This follows from parts (c), (e) and (f).

\smnoind
(h) If $r\in \BQ$ and $\bp\in \KPI$, then $|\bp(r)|_q^{1/n} \to 1$ when $\bp(r)\neq 0$.
Thus, \eqref{eqt:val-on-BQ} tells us that $\mu_{\BQ_q^{1/n}}(r) \to \mu_{\BQ_0}(r)$.
\end{prf}

\medskip

To simply the notation, if $K$ is  a minimal field, we may ignore $\mu_K$ and identify $K$ with its homeomorphic image in $\TBKM$, although this may occasionally cause ambiguity.

\medskip

As in the paragraph following Remark \ref{rem:neigh-TBK-M}, $\BQ_0$ is not closed in $\TBK$.
Thus, the closedness of $\BQ_0$ in $\TBKM$ (as established in Lemma \ref{lem:top-K-min}(a)) will prevent $\TBKM$ to be closed in $\TBK$.

\medskip

The following example gives a complete picture of $(\AZ)_0= \{\lambda\in \TBK: |\lambda| = 0\}$.

\medskip

\begin{eg}\label{eg:TBK_0}
For any $p\in \BP$, $\upsilon\in (0,1]$ and $\omega\in (0,\infty)$, we set $0_\BQ$, $0_\BR^\upsilon$, $0_p^\omega$ and $0_p^\infty$ to be the zeros of $\BQ_0$,  $\BR_\upsilon$,  $\BQ_p^\omega$ and  $\BZ(p)$, respectively.
By Lemma \ref{lem:top-K-min}(b), these zeros are all distinct elements in $\TBKM$.
We denote $L:=\{0_\BR^\upsilon: \upsilon\in(0,1]\}$, $S:=\{0_p^\omega: p\in \BP;\omega\in (0,\infty)\}$ and $D:= \{0_p^\infty: p\in \BP\}$.
Since $(K,0)\sim (K^0, 0)$ when $K\in \CK$ and $K^0$ is the smallest closed subfield of $K$, one knows $(\AZ)_0 = L\cup S\cup D \cup \{0_\BQ\}$.

On the other hand, it is well-known that there is a bijection from $\CK^{\min}$ to $\KM(\BZ_1)$ sending $K$ to $\lambda_K$, where $\lambda_K(m) = |m|_K$ ($m\in \BZ$).
It is easy to check that the induced map from $(\AZ)_0$ to $\KM(\BZ_1)$ (via Lemma \ref{lem:top-K-min}(b)) is a homeomorphism.
Thus, $(\AZ)_0$ can be identified with following compact subset of $\BR^2$:
$$
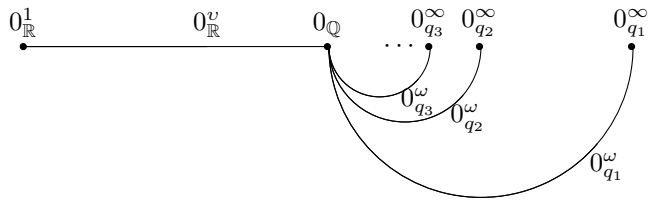
\begin{tikzpicture}[scale=2]

\draw (0,0) -- (1.2,0) node[above, black]{$0_\BR^\upsilon$};

\draw (0,0) -- (2,0);

\draw[fill=black] (-0.01,0) circle(.02) node[above, black]{$0_\BR^1$};

\draw[fill=black] (1.99,0) circle(.02) node[above, black]{$0_\BQ$};

\draw (2,0) arc(180:320:1/3) node[below, black]{$0_{q_3}^\omega$};

\draw (2,0) arc(180:360:1/3);

\draw[fill=black] (2+2/3-0.01,0) circle(.02) node[above, black]{$0_{q_3}^\infty$}  node[left, black]{$\cdots$};

\draw (2,0) arc(180:320:1/2) node[below, black]{$\ 0_{q_2}^\omega$};

\draw (2,0) arc(180:360:1/2);

\draw[fill=black] (2.99,0) circle(.02)  node[above, black]{$0_{q_2}^\infty$};

\draw (2,0) arc(180:320:1) node[below, black]{$\ \ 0_{q_1}^\omega$};

\draw (2,0) arc(180:360:1);

\draw[fill=black] (3.99,0) circle(.02)  node[above, black]{$0_{q_1}^\infty$};

\end{tikzpicture}$$
where $\BP = \{q_1, q_2, q_3, \cdots\}$ with $q_j$ being arranged in the increasing order.
One  can also obtain the above by using the discussion as in Lemma \ref{lem:top-K-min}.

Let us give a clearer description of the convergences of sequences in $S$ to $0_\BQ$.
By the definition, the collection of sets
$$\{0_\BQ\}\cup \big\{0_p^\omega\in S: 1-\epsilon < |k|_p^\omega\big\} \qquad (k\in \BN;\epsilon \in (0,1))$$
forms a neighborhood base for $0_\BQ$ in $\{0_\BQ\}\cup S$.
On the other hand, if we set $\BP_k:= \{p\in \BP: p\leq k\}$, then the collection of sets of the form
\begin{equation}\label{eqt:loc-neigh-basis-in-S}
\{0_\BQ\}\cup \big\{0_p^\omega\in S: \omega\in (0,\infty);p\in \BP\setminus \BP_k\big\}\cup \bigcup_{p\in \BP_k} \big\{0_p^\omega\in S: 1-\epsilon < p^{-k\omega}\big\}
\end{equation}
for arbitrary $k\in \BN$ and $\epsilon \in (0,1)$ also forms a neigborhood base.
%\begin{comment}
%[The following argument will not be in the publised version.]
In fact, the set in \eqref{eqt:loc-neigh-basis-in-S} is clearly contained in $\{0_\BQ\}\cup \big\{0_p^\omega\in S: 1-\epsilon < |k|_p^\omega\big\}$ (note that if $p^l$ divides $k$, then $l\leq k$).
Conversely, if we put $N_k:= \prod_{q\in \BP_k} q^k$, then $p$ divides $N_k$ if and only if $p\in \BP_k$.
Moreover, $|N_k|_p = p^{-k}$ when $p\in \BP_k$.
Thus, the set in \eqref{eqt:loc-neigh-basis-in-S} actually equals $\{0_\BQ\} \cup \big\{0_p^\omega\in S: 1-\epsilon < |N_k|_p^\omega\big\}$.
%\end{comment}

Hence, if we define $\omega_{k,\epsilon,p}:= \frac{\ln(1-\epsilon)}{-k\ln p}$, then $0_{p_i}^{\omega_i} \to 0_\BQ$ if and only if for any $k\in \BN$ and $\epsilon > 0$, there is $i_0\in \BN$ such that whenever $i\geq i_0$, either $p_i > k$ or $\omega_i < \omega_{k, \epsilon, p_i}$.
%Note also that if $\epsilon\in (0,1/2)$, then $\omega_{k, \epsilon, p} \leq 1$ for any $k$ and $p$.
\end{eg}

\medskip

Notice that the canonical inclusion $\BZ\subseteq \KPI$ defines a continuous map 
$$\kp:\TBK = \KM^0(\KPI)\to \KM^0(\BZ)\cong\KM(\BZ_1).$$ 
%Note that $\kp$ is also open since it will send open subsets of the form as in \eqref{eqt:neigh-TBK} to open subsets of $\KM(\BZ_1)$. 
On the other hand, by considering $\BZ$ as a quotient ring of $\KPI$ in the canonical way, one also has  a continuous map $\iota:\KM^0(\BZ)\to \KM^0(\KPI)$ with $\kp\circ \iota = \id_{\KM(\BZ_1)}$. 
This gives a ``fibration'' of $\TBK$ over $\KM(\BZ_1)$.  
More precisely, $\TBK$ can be decomposed as:
$$\TBK = {\bigcup}_{K\in \CK^{\min}} \BA_{\BZ_1}^{1,K}, $$
where $\BA_{\BZ_1}^{1,K}:= \kp^{-1}(\mu_K(0))$, i.e.\
$$\BA_{\BZ_1}^{1,K} = \bigcup \{\mu_L(t): K \text{ is the smallest closed subfield of } L; t\in L\}.$$
In particular, $\BA_{\BZ_1}^{1,\BQ_0}$ (respectively, $\BA_{\BZ_1}^{1,\BZ(p)}$) is the images of all elements in all fields with characteristic zero (respectively, characteristic $p\in \BP$), equipped with the trivial valuations.
Note also that $\BA_{\BZ_1}^{1,\BQ_0}$ and $\BA_{\BZ_1}^{1,\BZ(p)}$ ($p\in \BP$) are closed subsets of the compact space $(\AZU)_1$.
Furthermore,
$$\AZU=\BA_{\BZ_1}^{1,\BQ_0}\cup \bigcup_{p\in \BP} \BA_{\BZ_1}^{1,\BZ(p)} \cup \bigcup_{p\in \BP} \bigcup_{\omega\in (0,\infty)} \BA_{\BZ_1}^{1,\BQ_p^\omega}.$$

\medskip

The following proposition follows from Lemma \ref{lem:top-K-min} and its argument. 

\medskip

\begin{prop}
(a) Let $q\in \BP$.
Both $\bigcup_{\upsilon\in (0,1]} \BA_{\BZ_1}^{1,\BR_\upsilon}$ and $\bigcup_{\omega\in (0,\infty)} \BA_{\BZ_1}^{1,\BQ_q^\omega}$ are open subsets of $\TBK$ and the intersection of their closures is a subset of $\BA_{\BZ_1}^{1,\BQ_0}$.
Moreover, the closure of $\bigcup_{\upsilon\in (0,1]} \BA_{\BZ_1}^{1,\BR_\upsilon}$ is contained in $\BA_{\BZ_1}^{1,\BQ_0} \cup \bigcup_{\upsilon\in (0,1]} \BA_{\BZ_1}^{1,\BR_\upsilon}$.

\smnoind
(b) $\bigcup_{\upsilon\in (0,1]} \BA_{\BZ_1}^{1,\BR_\upsilon}\cong (0,1]\times \mathbb{H}$ as topological spaces.
\end{prop}
%\begin{comment}
[This proof will not appear in the published version.]

\begin{prf}
(a) This follows directly from parts (c), (d) and (h) of Lemma \ref{lem:top-K-min}.

\smnoind
(b) Let $\BC_\upsilon$ be the field $\BC$ equipped with the valuation $|c|_{\BC_\upsilon} = |c|_\E^\upsilon$ ($c\in \BC$).
It is well-known that the only non-trivial valuation field extension of $\BR_\upsilon$ is $\BC_\upsilon$.
Moreover, since $\kp(s) = 0_\BR^\upsilon$ for any $s\in \BC_\upsilon$, we know that the images of all such $\mu_{\BC_\upsilon}(\BC_\upsilon)$ with different $\upsilon$ are disjoint.
As in Example \ref{eg:non-inj}(b), there is a homeomorphism  $\varphi_\upsilon: \mathbb{H} \to \mu_{\BC_\upsilon}(\BC_\upsilon)$ such that $\mu_{\BC_\upsilon} = \varphi_\upsilon\circ Q_{\mathbb{H}}$.
This produces a bijection $\Psi: (0,1]\times \mathbb{H} \to \bigcup_{\upsilon\in (0,1]} \BA_{\BZ_1}^{1,\BR_\upsilon}$ with $\Psi(\upsilon, Q_{\mathbb{H}}(c)) = \mu_{\BC_\upsilon}(c)$.
It is not not to see that $\Psi$ is continuous and it follows from a similar argument as that of Lemma \ref{lem:top-K-min}(c) that $\Psi$ is a homeomorphism.
\end{prf}
%\end{comment}

\medskip

Let us end this section by giving the following clear and complete picture of $\TBKM$.
Since its argument is tedious and lengthy, we will not include it in the published version (please see the arXiv version of the article for details).

\medskip

\begin{thm}\label{thm:descr-K-min}
Let $q\in \BP$, $m\in \BZ$ and $n\in \BN$ with $m$ and $n$ being relatively prime.
Let $\Phi_q:\BZ\to \BZ(q)$ be the quotient map.
If $q$ does not divide $n$, we consider $[\frac{m}{n}]_{q}$ to be the element in $\BZ(q)$ satisfying $\Phi_{q}(n) [\frac{m}{n}]_{q} = \Phi_{q}(m)$.

\smnoind
(a) $\TBKM = \bigcup_{\upsilon\in (0,1]} \BR_\upsilon \cup \BQ_0 \cup \bigcup_{p\in \BP} \BZ(p) \cup \bigcup_{p\in \BP}\bigcup_{\omega\in (0,\infty)} \BQ_p^\omega$ is a fiber space over $(\AZ)_0$ with the fibers over $0_\BR^\upsilon$, $0_\BQ$, $0_p^\omega$ and $0_p^\infty$ being the topological spaces $\BR_1$, $\BQ_0$, $\BQ_p^1$ and $\BZ(p)$, respectively.
Moreover, the fiber topologies on $\bigcup_{\upsilon\in (0,1]} \BR_\upsilon$ and  $\bigcup_{\omega\in (0,\infty)} \BQ_q^\omega$ ($q\in \BP$)  are the product topologies.

\smnoind
(b) The subsets $\bigcup_{\upsilon\in (0,1]} \BR_\upsilon$,  $\bigcup_{\omega\in (0,\infty)} \BQ_q^\omega$ and $\BZ(q) \cup \bigcup_{\omega\in (0,\infty)} \BQ_q^\omega$ are open in $\TBKM$.

\smnoind
(c) The topologies on $\BQ_0$ and $\bigcup_{p\in \BP} \BZ(p)$ are discrete.

\smnoind
(d) Suppose that $\{(p_i, \omega_i)\}_{i\in \BN}$ is a sequence in $\BP\times (0,\infty)$ and $s_i\in \BQ_{p_i}^{\omega_i}$ ($i\in \BN$).
Then $s_i \to \Phi_q(n)$ inside $\TBKM$ if and only if $\omega_i\to \infty$ and there exists $i_0\in \BN$ such that for any $i\geq i_0$, one has $p_i = q$ as well as $s_i \in n + q\cdot\BZ_q$, where
$\BZ_q:=\{s\in \BQ_q: |s|_q \leq 1\}$.

\smnoind
(e) Suppose that  $\{(p_i, k_i)\}_{i\in \BN}$ is a sequences in $\BP\times \BZ$.
Then $\Phi_{p_i}(k_i) \to \frac{m}{n}\in \BQ_0$ inside $\TBKM$ if and only if $p_i\to \infty$ and there is $i_0\in \BN$ with $\Phi_{p_i}(k_i) = [\frac{m}{n}]_{p_i}$ whenever $i\geq i_0$.

\smnoind
(f) Suppose that $\{\upsilon_i\}_{i\in \BN}$ is a sequence in $(0,1]$ and $t_i\in \BR^{\upsilon_i}$ ($i\in \BN$).
Then $t_i \to \frac{m}{n}\in \BQ_0$ inside $\TBKM$ if and only if $\upsilon_i \to 0$ and $|t_i - \frac{m}{n}|_\E^{\upsilon_i} \to 0$.

\smnoind
(g) Suppose that $\{(p_i, \omega_i)\}_{i\in \BN}$ is a sequence in $\BP\times (0,\infty)$ and $s_i\in \BQ_{p_i}^{\omega_i}$ ($i\in \BN$).
Then $s_i\to \frac{m}{n}\in \BQ_0$ inside $\TBKM$ if and only if
$|k|_{p_i}^{\omega_i}\to 1$ for any $k\in \BN$ and $|s_i-\frac{m}{n}|_{p_i}^{\omega_i} \to 0$.
\end{thm}
%\begin{comment}
[This proof will not appear in the published version.]

\begin{prf}
(a) This part follows from the discussion above (in particular, parts (c) and (d) of Lemma \ref{lem:top-K-min}).

\smnoind
(b) This follows from parts (c), (d) and (e) of Lemma \ref{lem:top-K-min} as well as Example \ref{eg:TBK_0}.

\smnoind
(c) This part follows from Proposition \ref{prop:top-K-alpha}(c).

\smnoind
(d) $\Rightarrow)$.
Since $0_{p_i}^{\omega_i} = \kp(s_i) \to \kp(\Phi_q(n)) = 0_q^\infty$, Example \ref{eg:TBK_0} tells us that $p_i = q$ when $i$ is bigger than a certain fixed integer $i_1$ and $\omega_i\to \infty$.
Since $|s_i - n|_q^{\omega_i} \to 0$ (by considering $\bt -n \in \KPI$), one can find $i_0 \geq i_1$ such that $|s_i - n|_q < 1$ whenever $i\geq i_0$.
Thus, $s_i \in n + q\cdot\BZ_q$ for any $i\geq i_0$.

\noindent
$\Leftarrow)$.
We may assume that $i_0 =1$.
Notice, first of all that, by using the argument of Lemma \ref{lem:top-K-min}(e),
\begin{equation}\label{eqt:bp(n)}
|\bp(n)|_q^{\omega_i} \to |\bp(\Phi_q(n))|_{\BZ(q)} \qquad (\bp\in \KPI).
\end{equation}
Let $\bp\in \KPI$ and $t_i\in q\cdot \BZ_q$ with $s_i = n +t_i$.
One can find $N\in \BN$ and $\bq_1,...,\bq_N\in \KPI$ satisfying
$$\bp(n+t_i) = \bp(n) + t_i\bq_1(n) + \cdots + t_i^N\bq_N(n).$$
Since $|t_i|_q \leq 1/q$, we know from \eqref{eqt:bp(n)} that $|t_i^k\bq_k(n)|_q^{\omega_i} = |t_i|_q^{k\omega_i} |\bq_k(n)|_q^{\omega_i} \to 0$ ($k=1,...,N$), and hence
$$\left|\sum_{k=1}^N t_i^k\bq_k(n)\right|_q^{\omega_i}\ \leq\ \max_{k=1,...,N} |t_i^k\bq_k(n)|_q^{\omega_i}\ \to\ 0.$$
If $q$ divides $\bp(n)$, then \eqref{eqt:bp(n)} implies that
$$|\bp(n+t_i)|_q^{\omega_i}\ \leq\ \max \left\{|\bp(n)|_q^{\omega_i}, \left|\sum_{k=1}^N t_i^k\bq_k(n)\right|_q^{\omega_i}\right\}\ \to\ 0\ =\ |\bp(\Phi_q(n))|_{\BZ(q)}.$$
If $q$ does not divide $\bp(n)$, then \eqref{eqt:bp(n)} gives $|\bp(n)|_q^{\omega_i}\to 1$.
Thus, when $i$ is large enough, one knows that $|\bp(n)|_q^{\omega_i}\neq \left|\sum_{k=1}^N t_i^k\bq_k(n)\right|_q^{\omega_i}$ and
$$|\bp(n+t_i)|_q^{\omega_i}\ =\ \max \left\{|\bp(n)|_q^{\omega_i}, \left|\sum_{k=1}^N t_i^k\bq_k(n)\right|_q^{\omega_i}\right\}.$$
This shows that $|\bp(n+t_i)|_q^{\omega_i} \to 1 = |\bp(\Phi_q(n))|_{\BZ(q)}$.

\smnoind
(e) $\Rightarrow)$.
By considering $\kp$, we know that $p_i\to \infty$.
Thus, $p_i$ is not a prime factor of $n$ (and $[\frac{m}{n}]_{p_i}$ makes sense) when $i$ is large enough.
Since $|n \Phi_{p_i}(k_i) - m|_{\BZ(p_i)} \to 0$ (by considering $n \bt -m \in \KPI$), we know that $\Phi_{p_i}(n) \Phi_{p_i}(k_i) = \Phi_{p_i}(m)$ eventually.

\smnoind
$\Leftarrow)$.
Take any $\bp\in \KPI$.
If $\bp(\frac{m}{n}) = 0$, then $\big|\bp\big(\Phi_{p_i}(k_i)\big)\big|_{\BZ(p_i)} = 0$ for $i \geq i_0$.
If $\bp(\frac{m}{n}) = \frac{m'}{n'}$ for some $m'\in \BZ\setminus \{0\}$ and $n'\in \BN$ with $m'$ and $n'$ being relatively prime, there exists $i_1\geq i_0$ with $p_i$ not being a prime factor of $m'$ nor $n'$, and hence $\big|\bp\big(\Phi_{p_i}(k_i)\big)\big|_{\BZ(p_i)} = 1$, whenever $i\geq i_1$.
Consequently, $\mu_{\BZ(p_i)}\big(\Phi_{p_i}(k_i)\big) \to \mu_{\BQ_0}(\frac{m}{n})$.

\smnoind
(f) $\Rightarrow)$.
Observe that $\upsilon_i \to 0$ because $\kp\big(\mu_{\BR^{\upsilon_i}}(t_i)\big)\to \kp\big(\mu_{\BQ_0}(\frac{m}{n})\big)$ (c.f.\ the argument of Lemma \ref{lem:top-K-min}(c)).
Moreover, we have
$\big|t_i - \frac{m}{n}\big|_\E^{\upsilon_i}\ =\ n^{-\upsilon_i} |nt_i - m|_\E^{\upsilon_i}\ \to\ 0$ (as $n\bt -m \in \KPI$).

\noindent
$\Leftarrow)$.
Let $\bp\in \KPI$.
One can find $N\in \BN$, $M\in \BZ$ and $\bq\in \KPI$ such that
$$N\bp(t)\ =\ M + (t- m/n)\bq(t-m/n)$$ for any element $t$ in any field of characteristic $0$.
As $\big|t_i - \frac{m}{n}\big|_\E^{\upsilon_i} \to 0$, there is $i_0\in \BN$ such that $\big|t_i - \frac{m}{n}\big|_\E < 1$ whenever $i\geq i_0$.
Therefore, $\{|\bq(t_i- \frac{m}{n})|_\E\}_{i\in \BN}$ is bounded, and we have (since $\upsilon_i\to 0$)
\begin{equation}\label{eqt:upsilon:to-0}
N^{-\upsilon_i}\cdot |t_i - m/n|_\E^{\upsilon_i}\cdot |\bq(t_i- m/n)|_\E^{\upsilon_i}\ \to\ 0.
\end{equation}
When $\bp(\frac{m}{n}) = 0$, one has $\bp(t) = N^{-1}(t - \frac{m}{n})\bq(t- \frac{m}{n})$, and \eqref{eqt:upsilon:to-0} implies that $$|\bp(t_i)|_\E^{\upsilon_i}\ \to\ 0 = |\bp(m/n)|_{\BQ_0}.$$
When $\bp(\frac{m}{n}) \neq 0$, we have $M\neq 0$ and
\begin{eqnarray*}
\lefteqn{|M/N|_\E^{\upsilon_i} - |N^{-1}(t_i - m/n)\bq(s_i- m/n)|_\E^{\upsilon_i} \quad \leq \quad |\bp(t_i)|_\E^{\upsilon_i}}\\
& \qquad \qquad \qquad \qquad \qquad \leq &
|M/N|_\E^{\upsilon_i} + |N^{-1}(t_i - m/n)\bq(t_i- m/n)|_\E^{\upsilon_i}
\end{eqnarray*}
(note that $|\cdot|_\E^{\upsilon_i}$ is a norm).
As $|M/N|_\E^{\upsilon_i}\to 1$, \eqref{eqt:upsilon:to-0} tells us that
$$|\bp(t_i)|_\E^{\upsilon_i}\ \to\ 1 = |\bp(m/n)|_{\BQ_0}.$$

\smnoind
(g) $\Rightarrow)$.
As $0_{p_i}^{\omega_i} = \kp(s_i) \to 0_\BQ$, we know that $|k|_{p_i}^{\omega_i}\to 1$ ($k\in \BN\subseteq \KPI$).
Moreover, since $|ns_i - m|_{p_i}^{\omega_i}\to 0$ (consider $n\bt -m \in \KPI$), one concludes that $|s_i-\frac{m}{n}|_{p_i}^{\omega_i} = |n|_{p_i}^{-\omega_i}|ns_i - m|_{p_i}^{\omega_i} \to 0$.

\noindent
$\Leftarrow)$.
Let $\bp\in \KPI$.
One can find $N\in \BN$, $M\in \BZ$ and $\bq\in \KPI$ with $N\bp(t) = M + (t- \frac{m}{n})\bq(t-\frac{m}{n})$ for any element $t$ in any field of characteristic $0$.
If $\bq = \sum_{l=0}^{M_0} k_l \bt^l$ (where $k_0,...,k_{M_0}\in \BZ$), then
\begin{equation}\label{eqt:omega:to-0}
\big|(s_i- m/n)\bq(s_i- m/n)\big|_{p_i}^{\omega_i}\ \leq\ \max_{l=0,...,M_0} |k_l|_{p_i}^{\omega_i}\cdot\big|s_i- m/n\big|_{p_i}^{(l+1)\omega_i}\ \to\ 0.
\end{equation}
Suppose that $\bp(\frac{m}{n}) = 0$, i.e., $M=0$.
Then \eqref{eqt:omega:to-0} and the hypothesis gives
$$|\bp(s_i)|_{p_i}^{\omega_i}\ = \ |N|_{p_i}^{-\omega_i}
\cdot |(s_i - m/n)\bq(s_i-m/n)|_{p_i}^{\omega_i}\ \to\ 0\ =\  |\bp(m/n)|_{\BQ_0}.$$
Suppose that $\bp(\frac{m}{n}) \neq 0$, or equivalently,  $M\neq 0$.
Then $|\bp(s_i)|_{p_i}^{\omega_i} = |N|_{p_i}^{-\omega_i}
\cdot |M + (s_i - m/n)\bq(s_i-m/n)|_{p_i}^{\omega_i}$.
Since $|M|_{p_i}^{\omega_i} \to 1$, we know from \eqref{eqt:omega:to-0} that  when $i$ is large, $|M|_{p_i}^{\omega_i} \neq |(s_i - m/n)\bq(s_i-m/n)|_{p_i}^{\omega_i}$, which means that
$$|\bp(s_i)|_{p_i}^{\omega_i}\ =\ |N|_{p_i}^{-\omega_i}\cdot \max\big\{|M|_{p_i}^{\omega_i} , |(s_i - m/n)\bq(s_i-m/n)|_{p_i}^{\omega_i}\big\}$$
(as $|\cdot|_{p_i}^{\omega_i}$ is an ultrametric norm), which converges to
$1 = |\bp(m/n)|_{\BQ_0}$.
\end{prf}
%\end{comment}

\bigskip

\section{Berkovich spectra of elements in Banach Rings}

In this section, $R$ is a unital Banach ring with cardinality dominated by an infinite cardinal $\alpha$ and $a\in R$.
We set

\begin{eqnarray*}
\lefteqn{\sigma^\Ber_R(a) := \big\{\ti s \in \ti \BK
: \text{there exist }K\in \CK \text{ and }s\in K\text{ such that }}\\
&\qquad \qquad \qquad \qquad& (K,s)_\sim = \ti s,\ R\hat \otimes_\BZ K\neq (0) \text{ and } 1\otimes s - a\otimes 1\notin \KG(R\hat \otimes_\BZ K)\big\}
\end{eqnarray*}
(see Lemma \ref{lem:k-lin} for the meaning of $\KG(\cdot)$) and
\begin{eqnarray*}
\lefteqn{\sigma_R^\um(a) := \big\{\ti s \in \ti \BK
: \text{there exist }K\in \CK \text{ and }s\in K\text{ such that }}\\
&\qquad \qquad \qquad & (K,s)_\sim = \ti s,\ R\hat \otimes^\um_\BZ K\neq (0) \text{ and } 1\otimes s - a\otimes 1\notin \KG(R\hat \otimes_\BZ^\um K)\big\}
\end{eqnarray*}
By considering the canonical homomorphism from $R\hat \otimes_\BZ K$ to $R\hat \otimes^\um_\BZ K$, we know that $\sigma_R^\um(a)\subseteq \sigma^\Ber_R(a)$.
Note also that $R\hat \otimes_\BZ K$ could be zero and the non-zero requirement in the definitions above are necessary.

\medskip

\begin{rem}\label{rem:ring-sp-anoth-look}
(a) Let $\ti s\in \sigma^\Ber_R(a)$ and $(K,s)$ be as in the definition of $\sigma^\Ber_R(a)$ above.
By Remark \ref{rem:K-ti-s}(b), there is a surjective isometry $\Phi\in \CC\big(\CH(\lambda_{\ti s});K^s\big)$ with $\Phi\big(\varphi_{\lambda_{\ti s}}(\bt)\big) = s$.
This gives a map $\Psi\in \CC\big(R\hat \otimes_\BZ \CH(\lambda_{\ti s}); R\hat \otimes_\BZ K\big)$ satisfying $\Psi\big(1\otimes \varphi_{\lambda_{\ti s}}(\bt)\big) = 1\otimes s$.
Hence, we have
$$\sigma^\Ber_R(a) =  \big\{\lambda \in \TBK: 1\otimes \varphi_{\lambda}(\bt) - a\otimes 1 \text{ is not invertible in }R\hat \otimes_\BZ \CH(\lambda)\neq (0)\}$$
(note that the the non-zeroness of $R\hat \otimes_\BZ \CH(\lambda)$ is an assumption instead of a fact).

\smnoind
(b) Let $\ti s\in \sigma^\um_R(a)$ and $(K,s)$ be as in the definition of $\sigma^\um_R(a)$ above.
Then $R\hat \otimes^\um_\BZ K$ is a unital ultrametric $K$-Banach algebra, which implies that $K$ is non-Archimedean.
In particular, $\sigma_\BR^\um(r)$ is the empty set for any $r\in \BR$ (see Lemma \ref{lem:proj-ten}(a)).
Moreover, the consideration as in part (a) gives
$$\sigma_R^\um(a) = \big\{\lambda \in \AZU: 1\otimes \varphi_{\lambda}(\bt) - a\otimes 1 \text{ is not invertible in }R\hat \otimes_\BZ^\um \CH(\lambda)\neq (0)\big\}.$$
\end{rem}

\medskip

\begin{lem}\label{lem:upper-semi-cont}
(a) Let $T$ be a commutative unital Banach ring.
For any $x\in R\hat \otimes_\BZ T$ (respectively, $x\in R\hat \otimes_\BZ^\um T$), the map from $\KM(T)$ to $\RP$ that sends $\lambda$ to $\|(\id \otimes \varphi_\lambda)(x)\|_\wedge$ (respectively, $\|(\id \otimes \varphi_\lambda)(x)\|_\wedge^\um$) is upper semi-continuous.

\smnoind
(b) If we put $\TBK(R):= \{\lambda\in \TBK: R\hat \otimes_\BZ \CH(\lambda)\neq (0)\}$, then
$$\TBK(R) = \bigcup \left\{\BA_{\BZ_1}^{1,K}: K\in \CK^{\min}; R\hat\otimes_\BZ K\neq (0)\right\}$$
and is a closed subset of $\TBK$.
\end{lem}
\begin{prf}
(a) We will only consider $x\in R\hat \otimes_\BZ T$, as the other case follows from a similar argument.
Suppose that $\kappa\in (0,\infty)$ and $\lambda\in \KM(T)$ with $\|(\id \otimes \varphi_\lambda)(x)\|_\wedge < \kappa$.

Let us first assume that $x\in R\otimes_\BZ T$.
There exists $u_1,...,u_n\in R$ and $v_1,...,v_n\in T$ such that $x = \sum_{i=1}^n u_i\otimes v_i$ and $\sum_{i=1}^n \|u_i\|\lambda(v_i) < \kappa$.
When $\gamma$ is closed enough to $\lambda$, we have
$\|(\id \otimes \varphi_\gamma)(x)\|_\wedge \leq \sum_{i=1}^n \|u_i\|\gamma(v_i) < \kappa$, which gives the upper semi-continuity.

For a general element $x\in R\hat\otimes_\BZ T$, we consider $\epsilon := (\kappa - \|(\id \otimes \varphi_\lambda)(x)\|_\wedge)/4$.
There is $y \in R\otimes_\BZ T$ with $\|x-y\|_\wedge < \epsilon$.
As $\id\otimes \varphi_\gamma: R\hat\otimes_\BZ T \to R\hat\otimes_\BZ \CH(\gamma)$ is a contraction for any $\gamma\in \KM(T)$, one may use the above and a standard argument to show that $\|(\id \otimes \varphi_\gamma)(x)\|_\wedge < \kappa$ when $\gamma$ is closed enough to $\lambda$.

\smnoind
(b) Let $\lambda\in \TBK$ and $K:=\CH(\lambda)$.
Then $K^0:=\CH(\kp(\lambda))$ is the smallest closed subfield of $K$. We set $A:=R\hat\otimes_\BZ K^0$.
By Example \ref{eg:zero-ten-prod}(c) and Lemma \ref{lem:proj-ten}(b),
$$R\hat \otimes_\BZ K
\ =\ R\hat \otimes_\BZ (K^0\hat \otimes_{K^0} K)
\ =\ R\hat \otimes_\BZ (K^0\hat \otimes_\BZ K)
\ =\ A\hat \otimes_\BZ K
\ =\ A\hat \otimes_{K^0} K.$$
Obviously, if $R\hat\otimes_\BZ K\neq (0)$, then $A\neq (0)$.
Conversely, if $A\neq (0)$, then %$A$ is a unital $K^0$-Banach algebra (see Lemma \ref{lem:proj-ten}(a)) and
$A\hat \otimes_{K^0} K\neq (0)$ (one may employ Lemma \ref{lem:non-zero-ten-norm-sp}(a) if $K_0$ is non-Archimedean).
Therefore, $\lambda\in \TBK(R)$ if and only if $\kp(\lambda)\in \TBK(R)$.
This gives the first equality.

Moreover, for any $\gamma\in (\AZ)_0$, one has $\gamma\in \TBK(R)$ if and only if $\|(\id\otimes \varphi_\gamma)(1\otimes 1)\|_\wedge \geq 1/2$.
Thus, by applying part (a) to $T= \BZ_1$, one knows that $(\AZ)_0 \cap \TBK(R)$ is a closed subset of $(\AZ)_0\cong \KM(\BZ_1)$, and $\TBK(R) = \kp^{-1}\big((\AZ)_0\cap \TBK(R)\big)$ is a closed subset of $\TBK$.
\end{prf}

\medskip

The following shows that $\sigma_R^\Ber(a)$ is always compact, and it is non-empty if it satisfies a mild assumption.
%In other words, $\sigma_R^\Ber(a)$ is a non-empty compact set if and only if $R$ is standard.
The idea of the proof of the compactness essentially comes from \cite[Theorem 7.1.2]{Berk90}.

\medskip

\begin{thm}\label{thm:closed-ring-sp}
Let $R$ be a unital Banach ring and $a\in R$.

\smnoind
(a) $\sigma^\um_R(a)$ and $\sigma^\Ber_R(a)$ (when non-empty) are compact subsets of $\TBK$.

\smnoind
(b) $\sigma_R^\um(a)\neq \emptyset$ if and only if there exists $K \in \CK^{\um,\min}$ with $R\hat \otimes_\BZ^\um K \neq (0)$, or equivalently, there is a non-zero contractive additive map $\psi$ from $R$ to some $L\in \CK^{\um }$.

\smnoind
(c) If there exists a non-zero contractive additive map from $R$ to some $L\in \CK$, then $\sigma_R^\Ber(a)\neq \emptyset$.
\end{thm}
\begin{prf}
(a) Since the argument for the compactness of $\sigma^\Ber_R(a)$ and $\sigma^\um_R(a)$ are the same, we will only establish the compactness of $\sigma^\Ber_R(a)$.
Let $M\in \BN$ with $\|a\| < M$.
One may use Lemma \ref{lem:k-lin}(d) to show that $\sigma^\Ber_R(a) \subseteq (\AZ)_M$.
As the set $(\AZ)_M= \KM(\BZ_1\la M^{-1}\bt\ra)$ is compact (see \eqref{eqt:spec-Z-t}), we know that $(\AZ)_M\cap \TBK(R)$ is compact (by Lemma \ref{lem:upper-semi-cont}(b)), and it suffices to show that $\sigma^\Ber_R(a)$ is closed in $(\AZ)_M\cap \TBK(R)$ (notice that $\sigma_R^\Ber(a)\subseteq \TBK(R)$ by the definition).
Suppose that $\lambda\in (\AZ)_M\cap \TBK(R)\setminus \sigma^\Ber_R(a)$ and $z\in R\hat \otimes_\BZ \CH(\lambda)$ is the inverse of $1\otimes \varphi_\lambda(\bt) - a\otimes 1$.
There exist $u_1,...,u_m\in R$ and $\bp_1,...,\bp_m,\bq_1,...,\bq_m\in \KPI$ such that $\varphi_\lambda(\bq_i)\neq 0$ ($i=1,...,m$) and
$$\left\|z - \sum_{i=1}^m u_i \otimes \varphi_\lambda(\bp_i)\varphi_\lambda(\bq_i)^{-1} \right\|_\wedge
 < \|1\otimes \varphi_\lambda(\bt) - a\otimes 1\|^{-1}
$$
which gives
\begin{equation}\label{eqt:inv-app}
\left\|1\otimes 1 - \big(1\otimes \varphi_\lambda(\bt) - a\otimes 1\big)\sum_{i=1}^m u_i \otimes \varphi_\lambda(\bp_i)\varphi_\lambda(\bq_i)^{-1} \right\|_\wedge
 < 1
\end{equation}
Let us set $\bar \bp_i := \bp_i \bq_1 \cdots \bq_{i-1} \bq_{i+1} \cdots \bq_m$ ($i=1,...,m$) and $\bq := \bq_1 \cdots \bq_m$.
Consider
$$x:= 1\otimes \bq - (1 \otimes \bt  - a\otimes 1)\sum_{i=1}^m u_i \otimes \bar \bp_i \in R \hat \otimes_\BZ \BZ_1\la M^{-1}\bt\ra.$$
Then \eqref{eqt:inv-app} tells us that there exists $\kappa\in (0,\infty)$ with
$\|(\id\otimes \varphi_\lambda)(x)\| < \kappa < \lambda(\bq)$.
Thus, when $\gamma\in (\AZ)_M\cap \TBK(R)$ is near to $\lambda$, we have
$\|(\id\otimes \varphi_\gamma)(x)\| < \kappa$ (because of Lemma \ref{lem:upper-semi-cont}(a)) and $\kappa < \gamma(\bq)$.
This means that
$$\left\| 1\otimes 1 - \big(1\otimes \varphi_\gamma(\bt) - a\otimes 1\big) \sum_{i=1}^m u_i \otimes \varphi_\gamma(\bar \bp_i)\varphi_\gamma(\bq)^{-1}\right\| < 1$$
and hence $1\otimes \varphi_\gamma(\bt) - a\otimes 1$ is right invertible (by Lemma \ref{lem:k-lin}(d)).
Similarly, one can show that when $\gamma$ is close enough to $\lambda$, then
$\big(\sum_{i=1}^m u_i \otimes \varphi_\gamma(\bar \bp_i)\varphi_\gamma(\bq)^{-1}\big)\big(1\otimes \varphi_\gamma(\bt) - a\otimes 1\big)$ is invertible, which implies that $1\otimes \varphi_\gamma(\bt) - a\otimes 1$ is left invertible.

\smnoind
(b) First of all, suppose that such a $\psi:R\to L\in \CK^{\um }$ exists.
Let $L^0\subseteq L$ be the smallest closed subfield of $L$.
Since $\psi\otimes\id: R\hat \otimes_\BZ^\um L_0 \to L$ is non-zero, we see that $R\hat \otimes_\BZ^\um L_0 \neq (0)$.

Secondly, suppose that $A:=R\hat \otimes_\BZ^\um K\neq (0)$ for some $K\in \CK^{\um, \min}$.
Then $A$ is a unital ultrametric $K$-Banach algebra, and we denote by $\Phi\in \CC(R; A)$ the canonical map.
One may then apply \cite[Theorem 7.1.2(i)]{Berk90} to obtain $\gamma\in \BA_K^1$ such that $1\otimes \varphi_\gamma(\bt) - \Phi(a)\otimes 1 \notin \KG\big(A\hat\otimes^\um_{K} \CH(\gamma)\big)$.
Consider the canonical unital ring homomorphism
\begin{equation*}%\label{eqt:ring-to-field}
\Upsilon: \KPI \to K[\bt],
\end{equation*}
and define $\lambda:= \gamma(\Upsilon(\cdot))\in \AZ$.
One can find a map $\bar\Upsilon_\gamma\in \CC(\CH(\lambda);\CH(\gamma))$ such that $\varphi_\gamma(\bt) = \bar \Upsilon_\gamma (\varphi_\lambda(\bt))$.
By considering $\Phi\otimes \bar \Upsilon_\gamma: R\hat\otimes^\um_\BZ \CH(\lambda)\to A\hat\otimes^\um_{K} \CH(\gamma)$, we see that $\big(\CH(\lambda), \varphi_\lambda(\bt)\big)_\sim \in \sigma_R^\um(a)$.

Finally, suppose that $\sigma_R^\um(a)\neq \emptyset$.
By the definition, $R\hat \otimes_\BZ^\um L\neq (0)$ for some $L\in \CK^{\um }$, and one has $R\hat \otimes_\BZ^\um L^0\neq (0)$.
As $R\hat \otimes_\BZ^\um L^0$ is a non-zero ultrametric $L^0$-Banach space and all fields in $\CK^{\um,\min}$ are spherically complete, one can find a contractive $L^0$-linear map from $R\hat \otimes_\BZ^\um L^0$ to $L^0$ sending $1$ to $1$, and its composition with the canonical map from $R$ to $R\hat \otimes_\BZ^\um L^0$ is the required additive and contractive map.

\smnoind
(c) If $L$ is non-Archimedean, then part (b) implies that $\sigma_R^\um(a) \neq \emptyset$ and hence $\sigma_R^\Ber(a)\neq \emptyset$.

Suppose that $L$ is Archimedean.
We may assume that $L = \BC_\upsilon$ for some $\upsilon \in (0,1]$.
By the hypothesis, $A:= R\hat\otimes_\BZ L$ is a unital $L$-Banach algebra.
If we set 
$$\sigma^\BC_A(a\otimes 1):=\{s\in \BC: s-a\otimes 1 \notin \KG(A)\},$$
one may employ a standard argument to show that $\sigma^\BC_A(a\otimes 1)\neq \emptyset$, and the non-emptiness of $\sigma_R^\Ber(a)$ follows from 
$\mu_L(\sigma^\BC_A(a\otimes 1)) \subseteq \sigma_R^\Ber(a)$. 
\end{prf}

\medskip

The hypothesis of part (c) above is satisfied when $R\hat \otimes_\BZ K\neq (0)$ for some $K\in \CK^{\min}_\BR$, or when $R\hat \otimes_\BZ^\um K\neq (0)$ for some $K\in \CK^{\um,\min}$.

\medskip

\begin{thm}\label{thm:non-empty:ring-triv-val}
(a) If $R$ is a regular (see Proposition \ref{prop:reg-ring-ten-prod}) ultrametric unital Banach ring and $a\in R$, then $\sigma_R^\um(a)\neq \emptyset$.

\smnoind
(b) Let $R$ be a unital ring.
Then $\sigma_{R_0}^\um(a)\neq \emptyset$ for every $a\in R$, where $R_0$ is as in Section 2.
\end{thm}
\begin{prf}
(a) It follows from Theorem \ref{thm:closed-ring-sp}(b) and Proposition \ref{prop:reg-ring-ten-prod}(b).

\smnoind
(b) We regard $R$ as an abelian group under addition.
There is an increasing family $\{G_i\}_{i\in \KI}$ of finitely generated additive subgroups of $R$ such that $\bigcup_{i\in \KI} G_i = R$ and all $G_i$ contains the identity $1_R$ of the ring $R$.
For each $i\in \KI$, there exists $n_i\in \BN$ and a torsion subgroup of $G_i$ with  $G_i \cong \BZ^{n_i}\oplus H_i$.
Moreover, we let $\varphi_R:\BZ\to R$ be the map as in the proof of Proposition \ref{prop:reg-ring-ten-prod}.

We first consider the case when $\varphi_R(\BZ)$ is infinite.
Then $1_R\otimes 1_\BQ \neq (0)$ in $G_i\otimes_\BZ \BQ \cong \BQ^{n_i}$ ($i\in \KI$).
As $R\otimes_\BZ \BQ = \bigcup_{i\in \KI} G_i\otimes_\BZ \BQ$, we knows that $1_R\otimes 1_\BQ \neq (0)$ in $R\otimes_\BZ \BQ$.
Since $R_0\hat \otimes_\BZ^\um \BQ_0 \cong R\otimes_\BZ \BQ$ as abelian groups, the conclusion follows from Theorem \ref{thm:closed-ring-sp}(b).

Secondly, assume that $\varphi_R(\BZ)\cong \BZ/n\BZ$ for some $n\in \BN$.
Let $p\in \BP$ be a prime factor of $n$ and $n':=n/p$.
Fix any $i\in \KI$.
Notice that $1_R\in H_i$.
If $H_i$ is decomposed as $\BZ/n_1\BZ \oplus \cdots \oplus \BZ/n_M\BZ$ and $n'\cdot 1_R$ corresponding to $(\bar k_1,...,\bar k_M)\in \BZ/n_1\BZ \oplus \cdots \oplus \BZ/n_M\BZ$, then
$$k_j = \begin{cases}
n_j/p & \text{ if } p \text{ divides } n_j\\
0     & \text{ otherwise.}
\end{cases}$$
This shows that $n'\cdot 1_R\otimes 1_{\BZ(p)}$ is non-zero in $H_i\otimes_\BZ (\BZ/p\BZ)$.
Consequently, $n'\cdot 1_R\otimes 1_{\BZ(p)}$ is non-zero in $R\otimes_\BZ (\BZ/p\BZ)$, which is isomorphic to $R_0\hat \otimes_\BZ^\um \BZ(p)$ as abelian groups.
Now, the conclusion again follows from Theorem \ref{thm:closed-ring-sp}(b).
\end{prf}

\medskip

Note that part (b) above does not follows from part (a) since a unital ring with trivial norm needs not be regular (e.g. $\BZ/4\BZ$). 

\medskip

Suppose that $K$ is a field and $A$ is a unital $K$-algebra. 
For any $a\in A$, we set
\begin{equation}\label{eqt:defn-ord-sp}
\sigma_A^K(a)\ :=\ \{t\in K: a - t\notin \KG(A)\}. 
\end{equation}
If $L\in \CK$, $B$ is a unital $L$-Banach algebra and $b\in B$, then  $\mu_L(\sigma_B^L(b))\subseteq \sigma_B^\Ber(b)$.
In fact, by considering the canonical map from $B\hat{\otimes}_\BZ L$ to $B$, it is easy to see that $B\hat{\otimes}_\BZ L\neq (0)$ and $b\otimes 1 - 1\otimes r\notin \KG(B\hat{\otimes}_\BZ L)$ whenever $r\in \sigma_B^L(b)$. 
Similarly,  $\mu_L(\sigma_B^L(b))\subseteq \sigma_B^\um(b)$ when $B$ is ultrametric. 

\medskip

\begin{eg}\label{eg:non-empty-reg-ultra}
(a) Let $\Gamma$ be a discrete group and $K\in \CK$. 
	If $\ell^1(\Gamma;K)$ is the Banach ring of all absolutely summable maps from $\Gamma$ to $K$ (equipped with the $\ell^1$-norm) and $a\in \ell^1(\Gamma;K)$, then $\sigma^\um_{\ell^1(\Gamma;K)}(a)\neq \emptyset$. 
	Indeed, since $\ell^1(\Gamma;K)$ contained a base satisfying the condition as in \eqref{eqt:base} with $\kappa = 1$, the conclusion follows from Lemma \ref{lem:non-zero-ten-norm-sp}(b) and Theorem \ref{thm:closed-ring-sp}(b). 

\smnoind
(b) Let $A$ be a unital complex algebra and $a\in A$. 
Then $\mu_{\BC_0}(\sigma_A^\BC(a))\subseteq \sigma_{A_0}^\um(a)$ (see the second paragraph of Section 2 for the notation $A_0$ and $\BC_0$). 
Note that $\mu_{\BC_0}(\sigma_A^\BC(a))$ is a subset of $\BA_{\BZ_1}^{1,\BQ_0}$ instead of $\BA_{\BZ_1}^{1,\BR_1}$. 
\end{eg}

\medskip

In the remainder of this section, we consider $T$ to be a commutative unital Banach ring with its cardinality not exceeding $\alpha$.
The following lemma is more or less well-known.

\medskip

\begin{lem}\label{lem:M(A)}
(a) $\lambda\mapsto (\CH(\lambda), \varphi_\lambda)$ (see Proposition \ref{prop:mult-norm>field}(a)) sets up a bijection from $\KM(T)$ onto $\Lambda_T/\approx$, where
\begin{eqnarray*}
\lefteqn{\Lambda_T:=\big\{(K, \varphi)
: K\in \CK; \varphi\in \CC(T;K) \text{ such that }}\\
&\qquad \qquad \qquad & K \text{ is generated by }\varphi(T) \text{ as a complete valuation field}\big\}
\end{eqnarray*}
and $(K,\varphi)\approx (L,\psi)$ if there is a bijection $\theta\in \CC(K;L)$ satisfying $\psi = \theta\circ \varphi$.

\smnoind
(b) If $T$ is ultrametric and $K\in \CK$ with $\CC(T;K)\neq \emptyset$, then $K$ is non-Archimedean.
\end{lem}

\medskip

For every $u\in T$, we define $\Phi_u: \KM(T) \to \AZ$ by
$$\Phi_u(\lambda)(\bp) := \lambda (\bp(u))\qquad (\lambda\in \KM(T);\bp\in \KPI).$$
If we regard $\KM(T) = \Lambda_T/\approx$ through Lemma \ref{lem:M(A)}(a) and identify $\AZ = \ti\BK$ as in the above, we have
$\Phi_u\big((K,\varphi)_{\approx}\big)= (K,\varphi(u))_\sim$ where $(K,\varphi)\in \Lambda_T$.
The following result can be regarded as an analogue of \cite[Proposition 7.1.4(i)]{Berk90}.

\medskip

\begin{prop}\label{prop:Omega->sigma}
Let $T$ be a commutative unital Banach ring and $u\in T$.

\smnoind
(a) $\Phi_u$ is continuous and $\Phi_u(\KM(T)) = \sigma^\Ber_T(u)$.
Thus, $\sigma_T^\Ber(u)\neq \emptyset$.

\smnoind
(b) If $T$ is generated by $u$ as a unital Banach ring, then $\Phi_u$ is a homeomoprhism.

\smnoind
(c) If $T$ is ultrametric, then $\Phi_u(\KM(T)) = \sigma_T^\um(u)$.
In particular, $\sigma_T^\um(u) = \sigma^\Ber_T(u)$.
\end{prop}
\begin{prf}
(a) Clearly, $\Phi_u$ is a continuous map.
Suppose that $\lambda\in \KM(T)$ and $K:= \CH(\lambda)$.
We first show that $\Phi_u(\lambda) \in \sigma^\Ber_T(u)$.
In fact, the map $\varphi_\lambda\in \CC(T;K)$ induces a map $\bar \varphi_\lambda\in \CC(T\hat \otimes_\BZ K; K)$ and we have $T\hat \otimes_\BZ K\neq (0)$.
Since
$$\bar \varphi_\lambda\big(u\otimes 1 - 1\otimes \varphi_\lambda(u)\big) = 0,$$
we know that $1\otimes \varphi_\lambda(u)- u\otimes 1 \notin \KG(T\hat\otimes_\BZ K)$.
Thus, $\Phi_u(\lambda) = (K, \varphi_\lambda(u))_\sim \in \sigma_T^\Ber(u)$.

Conversely, suppose that $(L,s)_\sim\in \sigma^\Ber_T(u)$.
By \cite[Corollary 1.2.4]{Berk90}, there is $H\in \CK$ and
$\psi\in \CC\big(T\hat\otimes_\BZ L; H\big)$ such that $\psi\big(u\otimes 1 - 1\otimes s\big) = 0$.
Let us define $\phi\in \CC(T;H)$ and $\chi\in \CC(L; H)$, respectively, by
$$\phi(a) := \psi(a\otimes 1) \quad \text{and}\quad \chi(r) := \psi(1\otimes r) \qquad (a\in T;r\in L).$$
Consider $L'$ to be the closed subfield of $H$ generated by $\phi(T)$.
Then $(L',\phi)\in \Lambda_T$.
Since $\chi(s) = \phi(u)$, we know that $(L,s) \sim \big(H, \chi(s)\big) \sim (L', \phi(u))$, which gives the surjectivity of $\Phi_u$.

\smnoind
(b) Since $\{\bp(u):\bp\in \KPI\}$ is dense in $T$, it follows from the definition of $\Phi_u$ that it is injective (and hence homeomorphic).

\smnoind
(c) By Lemma \ref{lem:M(A)}(b), we know that if $(K,\varphi)\in \Lambda_T$, then $K\in \CK^{\um }$.
Now, the argument of the first statement of part (a) gives $\Phi_u(\KM(T)) = \sigma_T^\um(u)$.
\end{prf}

\medskip

\begin{rem}\label{rem:sub-ring-non-empty}
Let $R$ be a unital Banach ring, $a\in R$ and $T$ is the unital Banach subring of $R$ generated by $a$.
It is natural to ask whether $\partial \sigma_T^\Ber(a) \subseteq \sigma_R^\Ber(a)$ or $\partial \sigma_T^\um(a) \subseteq \sigma_R^\um(a)$, where the boundaries are taken in $\TBK$.
Notice that if these inclusions hold, then one may use Proposition \ref{prop:Omega->sigma} to conclude the non-emptiness of $\sigma_R^\Ber(a)$ and $\sigma_R^\um(a)$.
However, we will see in Example \ref{eg:spec-diff-rings} below that neither of these inclusions hold in general.
\end{rem}

\medskip

In the remainder of this section, we will consider the case of a unital Banach algebra $A$ over a complete valuation $\kk$.
As in the above, we let $\alpha$ be an infinite cardinal larger than the cardinality of $A$.
We will regard $\kk$ as a unital Banach subring of $A$ in the usual way.
Let us set
$$\CK^\kk:= \{K\in \CK: K \text{ is a complete valuation field extension of }\kk\}.$$
For any $a\in A$, we define
\begin{eqnarray*}
\lefteqn{\sigma^\Ber_{A,\kk}(a) := \big\{\ti s \in \ti \BK
: \text{there exist }K\in \CK^\kk \text{ and }s\in K \text{ such that}}\\
&\qquad & \qquad \qquad \qquad (K,s)_\sim = \ti s \text{ and } 1\otimes s - a\otimes 1 \notin \KG(A\hat \otimes_\kk K)\big\},
\end{eqnarray*}
(observe that $A\hat \otimes_\kk K$ is always non-zero because of Lemma \ref{lem:non-zero-ten-norm-sp}(a)) as well as
\begin{eqnarray*}
\lefteqn{\sigma^\um_{A,\kk}(a) := \big\{\ti s \in \ti \BK
: \text{there exist }K\in \CK^\kk \text{ and }s\in K \text{ such that}}\\
&\qquad & \qquad \qquad (K,s)_\sim = \ti s, A\hat \otimes_\kk^\um K\neq (0),  \text{ and } 1\otimes s - a\otimes 1 \notin \KG(A\hat \otimes_\kk^\um K)\big\}
\end{eqnarray*}
(note that as $A$ is not assumed to be ultrametric, the requirement $A\hat \otimes_\kk^\um K\neq (0)$ is necessary).
One always has $\sigma_{A,\kk}^\um (a)\subseteq \sigma_A^\um (a)$, thanks to the canonical ring homomorphism from $A\hat \otimes^\um_\BZ K$ to $A\hat \otimes^\um_\kk K$, as well as
$\sigma_{A,\kk}^\um (a)\subseteq \sigma^\Ber_{A,\kk}(a) \subseteq \sigma^\Ber_A(a)$.

\medskip

On the other hand, we consider $\Upsilon: \KPI\to \kk[\bt]$ to be the canonical ring homomorphism and $\ti \Upsilon: \BA_\kk^1\to \AZ$ to be the induced continuous map.
If $\kk$ is minimum, then $\ti \Upsilon$ is easily seen to be injective.
However, this map is non-injective in general; for example, when $\kk = \BC$.

\smnoind

For any $\gamma\in \BA_\kk^1$, we actually have $\gamma(r) = |r|_\kk$ ($r\in \kk$), and hence $\CH(\lambda)\in \CK^\kk$.
Let us extend the notation in \cite{Berk90} slightly as follows:
$$\Sigma_A(a) :=  \big\{\gamma \in \BA_\kk^1: 1\otimes \varphi_\gamma(\bt) - a\otimes 1 \notin \KG\big(A\hat \otimes_\kk \CH(\gamma)\big)\big\}$$
as well as
$$\Sigma^\um_A(a) :=  \big\{\gamma \in \BA_\kk^1: A\hat \otimes^\um_\kk \CH(\gamma)\neq (0); 1\otimes \varphi_\gamma(\bt) - a\otimes 1 \notin \KG\big(A\hat \otimes^\um_\kk \CH(\gamma)\big)\big\}.$$
When $A$ is ultrametric, $A\hat \otimes_\kk^\um K\neq (0)$ for any $K\in \CK^\kk$ and $\Sigma^\um_A(a)$ coincides with the spectrum $\Sigma_a$ as defined in  \cite{Berk90}.
On the other hand, $\Sigma^\um_A(a) = \emptyset$ when $\kk$ is Archimedean.

\medskip

\begin{rem}\label{rem:alg-sp-anoth-look}
(a) Suppose that $K\in \CK^\kk$ and $s\in K$.
The map $\bq\mapsto |\bq(s)|_K$ ($\bq\in \kk[t]$) clearly belongs to $\BA_\kk^1$.
This gives a map $\mu_K^\kk: K\to \BA_\kk^1$.
Notice that $\mu_\kk^\kk$ is injective, and that $\mu_K = \ti \Upsilon \circ \mu_K^\kk$.

\smnoind
(b) One has
$$\Sigma_A(a)= \bigcup_{K\in \CK^\kk} \mu_K^\kk\big(\sigma_{A\hat\otimes_\kk K}^K(a\otimes 1)\big)
\quad \text{ and } \quad
\Sigma_A^\um(a)= \bigcup_{K\in \CK^\kk} \mu_K^\kk\big(\sigma_{A\hat\otimes_\kk^\um K}^K(a\otimes 1)\big)$$
(see \eqref{eqt:defn-ord-sp} for the meaning of $\sigma_B^K(b)$). 
\end{rem}

\medskip

The argument for the following result is standard and is left to the reader.

\medskip

\begin{lem}\label{lem:rel-bet-Sigma-and-Ber}
$\ti \Upsilon(\Sigma_A(a)) = \sigma_{A,\kk}^\Ber(a)$ and $\ti \Upsilon(\Sigma^\um_A(a)) = \sigma_{A,\kk}^\um(a)$.
\end{lem}
%\begin{comment}
[This proof will not appear in the published version.]

\begin{prf}
Suppose that $\gamma\in \Sigma_A(a)$. 
By Lemma \ref{lem:k-lin}(a), $\varphi_\gamma|_\kk\in \CC(\kk; \CH(\gamma))$ is isometric, and $\CH(\gamma)\in \CK^\kk$ with  $1\otimes \varphi_\gamma(\bt) - a\otimes 1 
\notin \KG(A\hat\otimes_\kk \CH(\gamma))$. 
Hence, $\big(\CH(\gamma), \varphi_\gamma(\bt)\big)_\sim\in \sigma_{A,\kk}^\Ber(a)$. 

Conversely, let $K\in \CK^\kk$ and $s\in K$ such that $1\otimes s - a\otimes 1 \notin \KG(A\hat\otimes_\kk K)$. 
If we define $\gamma(\bp) := |\bp(s)|_K$ ($\bp\in \kk[\bt]$), then $\gamma\in \BA_\kk^1$. 
Since $\varphi_\gamma(\bp)\mapsto \bp(s)$ is a well-defined isometric unital ring homomorphism from $\varphi_\lambda(\kk[\bt])$ to $K$, it extends to an isometry $\Psi\in \CC\big(\CH(\gamma); K)$. 
By considering $\id\otimes \Psi: A\hat\otimes_\kk \CH(\gamma) \to A\hat\otimes_\kk K$, one knows that $\gamma\in \Sigma_A(a)$.
Now, as $(\CH(\gamma),\varphi_\gamma(\bt))\sim (K,s)$, we conclude that $\ti\Upsilon(\gamma) = (K,s)_\sim$. 

The argument for the second equality is similar.
\end{prf}
%\end{comment}

\medskip

Moreover, by Lemma \ref{lem:non-zero-ten-norm-sp}(b) and a similar argument as Theorem \ref{thm:closed-ring-sp}, we obtain parts (a) to (d) of the following proposition.
Note that part (e) follows from part (c), while part (f) follows from Remark \ref{rem:alg-sp-anoth-look}(b) as well as a standard algbraic argument.

\medskip

\begin{prop}\label{prop:cpt-non-empty-alg-spec}
Let $A$ and $B$ be $\kk$-Banach algebras with $A$ being unital, and let $a\in A$.
Denote by $B\oplus^{1} \kk$ the unital $\kk$-algebra $B\oplus \kk$ equipped with the norm $\|(b,\lambda)\| = \|b\|_B+|\lambda|_\kk$.

\smnoind
(a) $\Sigma_A(a)$ and $\Sigma_A^\um(a)$ are compact subsets of $\BA_\kk^1$.

\smnoind
(b) $\Sigma^\um_A(a)\neq \emptyset$ if and only if $\kk$ is non-Archimedean and there is a non-zero contractive additive map from $A$ to $\kk$.

\smnoind
(c) If there is a non-zero contractive additive map from $A$ to $\kk$, then $\Sigma_A(a) \neq \emptyset$.

\smnoind
(d) If $\kk$ is minimum, $\Sigma_A(a)\cong \sigma^\Ber_{A,\kk}(a) = \sigma^\Ber_A(a)$ and $\Sigma^\um_A(a)\cong \sigma_{A,\kk}^\um(a) = \sigma_A^\um(a)$.

\smnoind
(e) $\Sigma_{B\oplus^{1} \kk}(d)\neq \emptyset$ for any $d\in B\oplus^{1} \kk$.

\smnoind
(f) If $B$ is untial, then
$\Sigma_{B\oplus^{1} \kk}(b) \setminus \{\mu_\kk^\kk(0)\} = \Sigma_B(b) \setminus \{\mu_\kk^\kk(0)\}$, for any $b\in B$.
\end{prop}
%\begin{comment}
[The following proof will not appear in the published version.]

\begin{prf}
(f) Note that for $K \in \CK^\alpha_\kk$, the $K$-algebra $A\hat\otimes_\kk K$ is the $K$-algebra unitalization of $B\hat\otimes_\kk K$. 
It follows from standard algebraic argument that 
$\sigma^K_{A\hat\otimes_\kk K}(b\otimes 1) \cup\{0\} = \sigma^K_{B\hat\otimes_\kk K}(b\otimes 1) \cup\{0\}$. 
Since $\mu_K^\kk(0) = \mu_\kk^\kk(0)$, Remark \ref{rem:alg-sp-anoth-look}(b) gives the required conclusion. 
\end{prf}
%\end{comment}

\medskip

In particular, if either $\kk$ is Archimedean, or $A$ is ultrametric, or $A$ is finite dimensional, or $A$ is commutative, then $\Sigma_A(a) \neq \emptyset$.
Moreover, $\Sigma_A^\um(a) \neq \emptyset$ when $\sigma^\um_{A}(a)$ is non-empty.
This fact follows from Theorem \ref{thm:closed-ring-sp}(b), Lemmas \ref{lem:proj-ten}(a), \ref{lem:non-zero-ten-norm-sp}(b) and \ref{lem:rel-bet-Sigma-and-Ber} as well as Propositions \ref{prop:cpt-non-empty-alg-spec}(b).

\medskip

In the following, $D$ is a commutative unital $\kk$-Banach algebra.
Suppose that $\lambda\in \KM(D)$ and $\varphi_\lambda\in \CC(D;\CH(\lambda))$ is as in Proposition \ref{prop:mult-norm>field}(a).
By Lemma \ref{lem:k-lin}(a), $\varphi_\lambda|_\kk$ is an isometry and we have
$\lambda(t) = |t|_\kk$ ($t\in \kk$).
Consequently, for any $b\in D$, one can define a map $\Psi_b: \KM(D) \to \BA_\kk^1$ by
$$\Psi_b(\gamma)(\bp):= \gamma(\bp(b))\qquad (\gamma\in \KM_\kk(D); \bp\in \kk[\bt]).$$
Clearly,
$\Phi_b = \ti \Upsilon \circ \Psi_b$.

\medskip

Using the argument of Proposition \ref{prop:Omega->sigma}, we obtain parts (a), (b) and (c) of the following result.
Notice also that the first statement of part (c) is precisely \cite[Proposition 7.1.4(a)]{Berk90}.
On the other hand, part (d) of this result follows from part (a), Lemma \ref{lem:rel-bet-Sigma-and-Ber} and Proposition \ref{prop:Omega->sigma}(a).
Similarly, part (e) follows from part (c), Lemma \ref{lem:rel-bet-Sigma-and-Ber} and \ref{prop:Omega->sigma}(c).
Furthermore, part (f) follows from the Mazur theorem, which states that if a commutative $\BR_1$-Banach algebra is also a division ring, then it is either $\BR_1$ or $\BC_1$.

\medskip

\begin{prop}\label{prop:Omega-alg->Sigma}
Let $D$ be a commutative unital $\kk$-Banach algebra and $b\in D$.

\smnoind
(a) $\Psi_b$ is a continuous surjection from $\KM(D)$ onto $\Sigma_D(b)$.

\smnoind
(b) If $D$ is generated by $b$ as a unital $\kk$-Banach algebra, $\Psi_b$ is a homeomoprhism.

\smnoind
(c) If $D$ is ultrametric, then $\Psi_b(\KM(D)) = \Sigma^\um_D(b)$.
In particular, $\Sigma^\um_D(b) = \Sigma_D(b)$.

\smnoind
(d) $\sigma^\Ber_{D,\kk}(b) = \sigma^\Ber_D(b)$.

\smnoind
(e) If $D$ is ultrametric, then $\sigma_{D,\kk}^\um (b) = \sigma^\um_D(b)$.

\smnoind
(f) If $\kk = \BR_1$, then $\Lambda_D= \CC(D;\BC_1)$.
\end{prop}

\medskip

Let us consider the cases when $\kk = \BR_1$ or $\BC_1$ and compare the above with the usual spectrum. 

\medskip

\begin{prop}\label{prop:spec-elem-in-real-alg}
Suppose that $B$ is a unital $\BR_1$-Banach algebra and $b\in B$. 

\smnoind
(a) $\sigma^\Ber_B(b) = \mu_{\BC_1}\Big(\sigma_{B\otimes_{\BR} {\BC}}^{\BC}(b\otimes 1)\Big)$.

\smnoind
(b) If $B$ is a unital $\BC_1$-Banach algebra, then $\sigma^\Ber_B(b) = \sigma^\Ber_{B,{\BC_1}}(b)= \mu_{\BC_1}\big(\sigma_B^{\BC}(b)\big)$.
\end{prop}
\begin{prf}
(a) This follows from Proposition \ref{prop:cpt-non-empty-alg-spec}(d), Remark \ref{rem:alg-sp-anoth-look}(b) as well as the facts that $\CK^{\BR_1} = \{{\BR_1},{\BC_1}\}$ and $B\hat \otimes_{\BR_1} \BC_1 = B\otimes_\BR \BC$ as complex algebras.

\smnoind
(b) The first equality follows from part (a), Example \ref{eg:non-inj}(b) and \cite[Proposition 3.3]{Ing64}.
The second equality follows from Remark \ref{rem:alg-sp-anoth-look}(b) and the fact that $\CK^{\BC_1} = \{{\BC_1}\}$.
\end{prf}

\medskip

Part (b) of the above tells us that if $B$ is a complex Banach algebra, then $\sigma^\Ber_B(b)$ is the ``folding up'' of the usual spectrum $\sigma_B^{\BC}(b)$ alone the real axis (see Example \ref{eg:non-inj}(b).

\medskip

Finally, we give a negative answer to the question raised in Remark \ref{rem:sub-ring-non-empty}.  

\medskip

\begin{eg}\label{eg:spec-diff-rings}
We use the notation as in Example \ref{eg:TBK_0}.

\smnoind	
(a) By Proposition \ref{prop:spec-elem-in-real-alg}, $\sigma_{\BR_1}^\Ber(0)  = \sigma_{\BC_1}^\Ber(0) = \{0_\BR^1\}$.
As the subring of $\BR_1$ generated by $0$ is $\BZ_1$, Proposition \ref{prop:Omega->sigma}(b) and the definition of $\Phi_0$ (as in Proposition \ref{prop:Omega->sigma} for $u=0$) tell us that  $\sigma_{\BZ_1}^\Ber(0) = (\AZ)_0$.
Consider $\omega\in (0,\infty)$ and $p\in \BP$.
Then $\mu_{\BQ_p^\omega}(p^n) \to 0_p^\omega$ when $n\to \infty$ and hence $0_p^\omega \in \partial \sigma_{\BZ_1}^\Ber(0) \setminus \sigma_{\BR_1}^\Ber(0)$ (the boundary is taken in $\TBK$).

\smnoind
(b) For any $K\in \CK^{\BQ_0}$, one has $\BQ_0\hat \otimes_{\BQ_0}^\um K = K$ and we know from Remark \ref{rem:ring-sp-anoth-look}(b) that $\sigma_{\BQ_0}^\um(0) = \{0_\BQ\}$.
The unital Banach subring of $\BQ_0$ generated by $0$ is $\BZ_0$.
Again Proposition \ref{prop:Omega->sigma}(b)\&(c) and the definition of $\Phi_0$ tells
us that $\sigma_{\BZ_0}^\um(0) = (\AZ)_0\cap \AZU$.
The same argument as part (a) implies that $0_p^\omega \in \partial \sigma_{\BZ_0}^\um(0) \setminus \sigma_{\BQ_0}^\um(0)$.
\end{eg}

\medskip

Nevertheless, we do not know whether $\partial \sigma_{B,\kk}^\Ber(a) \subseteq \sigma_{A,\kk}^\Ber(a)$, when $A$ is a unital $\kk$-Banach algebra, $a\in A$ and $B$ is the unital $\kk$-Banach subalgebra of $A$ generated by $a$.
Note that a positive answer to this question will give the non-emptiness of $\sigma_{A,\kk}^\Ber(a)$ (because of Propositions \ref{prop:Omega->sigma}(a) and  \ref{prop:Omega-alg->Sigma}(d)).

\medskip

\section{Existence of zeros of the Fredholm determinants}

\medskip

In the following, $\kk$ is a non-Archimedean complete valuation field and $\KJ$ is an infinite set.
Suppose that $E := c_0(\KJ;\kk)$ with $\{e_j\}_{j\in \KJ}$ being its canonical base.
We set $\alpha$ to be the cardinality of $\CL(E)$.
For any $v\in \CL(E)$, if $(\eta_{ij})_{i,j\in \KJ}$ is the matrix representing $v$ under the base $\{e_j\}_{j\in \KJ}$, then, as in the paragraph preceeding \cite[Proposition 4]{Serre62}, one has
\begin{equation}\label{eqt:norm-sup}
\|v\|_{\CL(E)}\ =\ \sup_{i,j\in \KJ} |\eta_{ij}|_\kk.
\end{equation}

\medskip

Let $u\in \CL(E)$ be a completely continuous operator and $\det (1-\bt u)\in \kk[[\bt]]$ be the Fredholm determinant as defined \cite[p.75]{Serre62}.
The aim of this section is to give an equivalent condition for the existence of a zero of $\det (1-\bt u)$ in a complete valuation field extension of $\kk$.

\medskip

We first start with the following well-known lemma.
Note that part (a) is precisely \cite[Proposition 6]{Serre62}, while part (b) follows from part (a) as well as \eqref{eqt:norm-sup}.

\medskip

\begin{lem}\label{lem:c_0-tensor-K}
Let $K\in \CK^\kk$.

\smnoind
(a) There is a natural isometric isomorphism from $E\hat\otimes_\kk^\um K$ to $c_0(\KJ;K)$ sending $\{e_j\otimes 1\}_{j\in \KJ}$ to the canonical base of $c_0(\KJ;K)$.

\smnoind
(b) There is a canonical isometric $\kk$-linear map $\kappa_{E,K}: \CL(E)\to \CL(E\hat\otimes_\kk^\um K)$ such that $\kappa_{E,K}(b)(x\otimes s) = b(x)s$ ($b\in \CL(E); x\in E; s\in K$).
\end{lem}

\medskip

As usual, we denote $\kappa_{E,K}(b)$ by $b\otimes \id$.
Observe that part (b) above and Lemma \ref{lem:univ-prop-proj-ten}(c) produce a contractive unital ring homomorphism $\hat \kappa_{E,K}: \CL(E)\hat\otimes_\kk^\um K \to \CL(E\hat\otimes_\kk^\um K)$ with
\begin{equation*}%\label{eqt:defn-hat-kappa}
\hat \kappa_{E,K}(b\otimes t) = (b\otimes \id)\cdot t \qquad (b\in \CL(E);t\in K)).
\end{equation*}
Moreover, we use $A_u$ to denote the unital $\kk$-Banach subalgebra of $\CL(E)$ generated by $u$.

\medskip

Our next lemma is more or less known, but we give its argument here for the sake of completeness.

\medskip

\begin{lem}\label{lem:0-in-sp}
	(a) $(\kk,0)_\sim\in \sigma^\Ber_{\CL(E), \kk}(u)$.
	
	\smnoind
	(b) If $(K,0)_\sim\in \sigma^\Ber_{\CL(E), \kk}(u)$, then
	$(K,0)_\sim = (\kk,0)_\sim$.

\smnoind
(c) If $s\in \kk$ such that $1-su \in \KG(\CL(E))$, then $(1-su)^{-1}\in A_u$.
\end{lem}
\begin{prf}
(a) It follows from the Corollary of \cite[Proposition 4]{Serre62} (which implies that the identity map is not completely continuous, and hence no completely continuous operator is invertible).

\smnoind
(b) This part follows from the fact that $K\in \CK^\kk$.

\smnoind
(c) Let $P(\bt,u)\in \CL(E)[[\bt]]$ be the Fredholm resolvent as in  \cite[p.78]{Serre62}.
Then
\begin{equation}\label{eqt:Fred-det}
\det (1-\bt u) = P(\bt,u)(1-\bt u)
\end{equation}
and for any $r\in \kk$, one knows that
$1-ru$ is invertible if and only if $\det (1-ru) \neq 0$ (see \cite[Proposition 11]{Serre62}).
Moreover, if we consider the expansion $\det (1-\bt u) = \sum_{i=0}^\infty c_i\bt^i$ and $P(\bt,u) = \sum_{i=0}^\infty x_i\bt^i$, where $c_0,c_1,c_2, ...\in \kk$ and $x_0,x_1,x_2,...\in \CL(E)$, then as in \cite[p.78]{Serre62},
\begin{equation}\label{eqt:Fred-inver}
x_0 = 1 \quad \text{and} \quad x_i = c_i + ux_{i-1} \qquad (i\in \BN).
\end{equation}
Thus, $x_i\in A_u$, for every $i\in \BN$, and $(1-su)^{-1} = \det (1-su)^{-1}P(s,u)\in A_u$.
\end{prf}

\medskip

\begin{prop}\label{prop:inv-in-A}
Let $E$, $u$ and $A_u$ be as in the above.
Then $\sigma^\um_{\CL(E),\kk}(u) = \sigma^\Ber_{\CL(E)}(u) = \sigma^\Ber_{A_u}(u) = \Phi_u(\KM(A_u))$ and they coincide with
\begin{eqnarray*}
\lefteqn{\sigma^\um(u) := \big\{\ti s \in \ti \BK: \text{there exist }K\in\CK^\kk \text{ and }s\in K\text{ such that } }\\
&\qquad \qquad \qquad & (K,s)_\sim = \ti s \text{ and } s - u\otimes \id \notin \KG\big(\CL(E\hat \otimes_\kk^\um K)\big)\big\}.
\end{eqnarray*}
\end{prop}
\begin{prf}
As in  Lemma \ref{lem:rel-bet-Sigma-and-Ber}, one has
\begin{equation}\label{eqt:anoth-form-sigma}
\sigma^\um(u) = \ti \Upsilon\left(\big\{\gamma \in \BA_\kk^1: \varphi_\gamma(\bt) - u\otimes \id \notin \KG\big(\CL(E\hat \otimes_\kk^\um \CH(\gamma))\big)\big\}\right).
\end{equation}
If $\gamma\in \BA_\kk^1$ satisfying $1\otimes \varphi_\gamma(\bt) - u\otimes 1 \in \KG\big(\CL(E)\hat\otimes_\kk^\um \CH(\gamma)\big)$, then $\varphi_\gamma(\bt) - u\otimes \id = \hat \kappa_{E,\CH(\gamma)}\big(1\otimes \varphi_\gamma(\bt) - u\otimes 1\big)$ is invertible in $\CL(E\hat \otimes_\kk^\um \CH(\gamma))$.
Thus, by Lemma \ref{lem:rel-bet-Sigma-and-Ber}, we know that $\sigma^\um(u)\subseteq \sigma^\um_{\CL(E),\kk}(u)$.

Obviously, $\sigma^\um_{\CL(E),\kk}(u)\subseteq \sigma^\Ber_{\CL(E)}(u) \subseteq \sigma^\Ber_{A_u}(u)$, and $\sigma^\Ber_{A_u}(u) = \Phi_u(\KM(A_u))$ because of Proposition \ref{prop:Omega->sigma}(a).
By Relation \eqref{eqt:anoth-form-sigma}, Proposition \ref{prop:Omega-alg->Sigma}(d) and Lemma \ref{lem:rel-bet-Sigma-and-Ber}, in order to verify $\sigma^\Ber_{A_u}(u)\subseteq \sigma^\um(u)$, it suffices to show that $\varphi_\gamma(\bt) - u\otimes \id \notin \KG\big(\CL(E\hat \otimes_\kk^\um \CH(\gamma))\big)$ whenever $\gamma\in \Sigma_{A_u}(u)$.

Suppose on the contrary that there exists $\lambda\in \Sigma_{A_u}(u)$ with $\varphi_\lambda(\bt) - u\otimes \id \in \KG\big(\CL(E\hat \otimes_\kk^\um \CH(\lambda))\big)$.
Since $u\otimes \id$ is completely continuous, Lemma \ref{lem:0-in-sp}(a) tells us that $r:= \varphi_\lambda(\bt)^{-1}$ makes sense.
As in Relation \eqref{eqt:Fred-det}, one has
$$(1 - r(u\otimes \id))^{-1} = \det (1-r(u\otimes \id))^{-1}P(r,u\otimes \id).$$
By Lemma \ref{lem:c_0-tensor-K}(a) and \cite[Proposition 7(a)]{Serre62}, the two Fredholm determinants $\det ( 1- \bt u)$ and $\det(1 - \bt (u\otimes \id))$  coincide.
Thus, if $\{c_i\}_{i\in \BN_0}$ and $\{x_i\}_{i\in \BN_0}$ are the sequences as in the argument of Lemma \ref{lem:0-in-sp}(c) (in particular, $x_i\in A_u$) and $P(\bt ,u\otimes \id)$ is expressed as $\sum_{i=0}^\infty y_i\bt^i$, where $y_i\in \CL(E\hat \otimes_\kk^\um \CH(\lambda))$, then one has
$$y_0 = 1 \quad \text{and} \quad y_i = c_i + (u\otimes \id)y_{i-1}\qquad (i\in \BN),$$
which gives $y_i = \kappa_{E,\CH(\lambda)}(x_i)$ ($i\in \BN$).
Since $P(r,u\otimes \id)$ converges, we know that $\|x_i\| |r|^i = \|y_i\| |r|^i \to 0$ as $i\to \infty$ (note that $\kappa_{E,\CH(\lambda)}$ is isometric).
Consequently, the element $\hat P(r,u\otimes 1) := \sum_{i=0}^\infty x_i\otimes r^i$ exists in $A_u\hat\otimes_\kk^\um \CH(\lambda)$ (as $\|\cdot\|_\wedge^\um$ is a subcross norm).
Moreover, as $1 - r(u\otimes \id)\in \KG\big(\CL(E\hat \otimes_\kk^\um \CH(\lambda))\big)$, we have
$$\sum_{i=0}^\infty c_ir^i\ = \ \det(1-r(u\otimes \id))\ \neq\ 0.$$
Now, we set
$b:= (\sum_{i=0}^\infty c_ir^i)^{-1} \hat P(r,u\otimes 1)$.
It is not hard to check, with the help of Equalities \eqref{eqt:Fred-inver}, that
$$\left(\sum_{i=0}^\infty c_ir^i\right)b(1-u\otimes r) = \left(\sum_{i=0}^\infty x_i\otimes r^i\right)(1-u\otimes r) = \sum_{i=0}^\infty c_ir^i = (1-u\otimes r)\left(\sum_{i=0}^\infty c_ir^i\right)b,$$
which means that $1-u\otimes r$ is both left and right invertible.
This gives the contradiction that $1\otimes \varphi_\lambda(\bt) - u\otimes 1\in \KG\big(A_u\hat \otimes_\kk^\um \CH(\lambda)\big)$.
\end{prf}

\medskip

One can use Equality \eqref{eqt:norm-sup} to show that if $\lambda\in \BA_\kk^1$ such that $\CH(\lambda)$, when considered as a $\kk$-Banach space, has an orthogonal base, then the canonical map from $A_u\hat\otimes^\um_\kk\CH(\lambda)$ to $\CL(E\hat \otimes^\um_\kk \CH(\lambda))$ is actually isometric, and the argument of the above can be simplified.
This applies, in particular, to the case when $\kk$ is separable and locally compact, becasue of \cite[Theorem 50.8]{Schik84}.

\medskip

Proposition \ref{prop:inv-in-A} produces a relationship between non-zero elements in the spectrum of $u$ and zeros of the Fredholm determinant associated with $u$, which, although not extremely surprising, gives a way to determine whether a zero of the Fredholm determinat exists.

\medskip

\begin{thm}\label{thm:zero-Fred-det}
Let $E$ and $u$ be as in the above.
Then
$$\left\{(K,s^{-1})_\sim\in \ti \BK: K\in \CK_\kk; s\in K; \det(1-su) = 0\right\}
\ =\ \sigma^\Ber_{\CL(E)}(u)\setminus \left\{(\kk, 0)_\sim\right\}.$$
\end{thm}
\begin{prf}
By Proposition  \ref{prop:inv-in-A} and Lemma \ref{lem:0-in-sp}(b), one has $(K,r)_\sim \in \sigma^\Ber_{\CL(E)}(u) \setminus \left\{(\kk, 0)_\sim\right\}$ if and only if $r\neq 0$ and $r - u\otimes \id \notin \KG\big(\CL(E\hat \otimes_\kk^\um K)\big)$.
Moreover, \cite[Proposition 11]{Serre62} tells us that the later is equivalent to $\det(1-r^{-1}u) = \det(1-r^{-1}(u\otimes \id)) = 0$ (see also the argument of Proposition  \ref{prop:inv-in-A}).
The equivalence now follows from the fact that $0$ is never a zero of $\det(1-\bt u)$ in any complete valuation field.
\end{prf}

\medskip

By \cite[Theorem 1.3.1]{Berk90} and Proposition \ref{prop:Omega->sigma}(a), we know that
$$\max \left\{|s|_K: (K,s)_\sim\in \sigma_{A_u}^\Ber(u)\right\}\ =\ \lim_n \|u^n\|^{\frac{1}{n}}.$$
This, together with Proposition \ref{prop:inv-in-A} and Theorem \ref{thm:zero-Fred-det}, gives the following corollary.

\medskip

\begin{cor}\label{cor:lower-bdd-zero}
Let $E$ and $u$ be as in the above. 

\smnoind
(a) If $K\in \CK_\kk$ and $s\in K$ satisfying $\det(1-su) = 0$, then $|s|_K \geq \lim_n \|u^n\|^{-1/n}$. 

\smnoind
(b) If $\lim_n \|u^n\|^{1/n} = 0$, there does not exist $K\in \CK_\kk$ and $s\in K$ with $\det(1-su) = 0$. 
\end{cor}

\medskip

Let us end this article with the following concrete example of a completely continuous operator whose Fredholm determinant has no zero in any complete valuation field extension.

\medskip

\begin{eg}\label{eg:Fredh-det-no-zero}
Let $p\in \BP$ and $E:= c_0(\BN; \BQ_p^1)$.
Define $u\in \CL(E)$ by
$$u\big((s_i)_{i\in \BN}\big) := (p^i\cdot s_{i+1})_{i\in \BN} \quad ((s_i)_{i\in \BN}\in E).$$
Clearly, $u$ is completely continuous.
It follows from Equality \eqref{eqt:norm-sup} that $\|u^n\| \leq p^{-\frac{n(n+1)}{2}}$ ($n\in \BN$).
Thus, Corollary \ref{cor:lower-bdd-zero}(b) tells us that there is no zero for the Fredholm determinant of $u$ in any complete valuation field extension of $\BQ_p^1$.
\end{eg}

%\medskip

%Moreover, $u\otimes\id \in \CL(E\hat\otimes_k^\um K)$ is again completely continuous.
%By Lemma \ref{lem:c_0-tensor-K}, it is easy to see that $\det (1-t(u\otimes \id)) = \det (1-tu)$.

\end{document}